\documentclass[english]{article}
\usepackage[T1]{fontenc}		
\usepackage[latin9]{inputenc}	
\usepackage{textcomp}	
\usepackage{amsmath}	
\usepackage{graphicx}   
\usepackage{epstopdf}   
\usepackage{amssymb}	
\usepackage[font=small,labelfont=bf,margin=24pt]{caption}	
\usepackage{url}		
\usepackage{esint}		
\usepackage{babel}		
\usepackage[usenames,dvipsnames]{pstricks}	
\usepackage{epsfig}		
\usepackage{pst-grad} 	
\usepackage{pst-plot} 	
\usepackage{mathtools} 	
\usepackage{geometry} 	

\usepackage{tikz} 						
\usetikzlibrary{decorations.pathreplacing}	
\usetikzlibrary{arrows}					

\usepackage{algorithm}
\usepackage{algpseudocode}

\usepackage{xpatch}

\graphicspath{ {img/}{phd_images/}{charts/} }

\algrenewcommand\alglinenumber[1]{{#1:\quad}}
\makeatletter
\xpatchcmd{\algorithmic}{\itemsep\z@}{\itemsep=1ex plus2pt}{}{}
\makeatother

\algnewcommand\algorithmicinput{\textbf{Input:}}
\algnewcommand\Input{\item[\algorithmicinput]} 
\algnewcommand\algorithmicoutput{\textbf{Output:}}
\algnewcommand\Output{\item[\algorithmicoutput]}

\newcommand{\repr}[4][x]{\ensuremath{{\bar{#1}}_{#4}^{#2,[#3]}}}
\newcommand{\epsi}{\ensuremath{\varepsilon}}

\begingroup
	 \pgfqkeys{/tikz}{active/.style={opacity=1, fill opacity = 0.2, text opacity=1}}			
	 \pgfqkeys{/tikz}{inactive/.style={opacity=0.1, fill opacity = 0.05, text opacity=0.2}}	
\endgroup

\newcommand{\reprnode}[4][]{
	\begin{scope}[{#1}]
		\newcommand\drcnn{\the\numexpr\reprn+1\relax}	
		\newcommand\drci{#2}		
		\newcommand\drck{#3}		
		\node at (-\drci / \reprp, \drcnn / \reprp * -\drck / \drcnn) ({#4}) {} ;
	\end{scope}
} 

\newcommand{\drawreprcolumn}[3][]{
	\begin{scope}[{#1}]
		\newcommand\drcnn{\the\numexpr\reprn+1\relax}	
		\newcommand\drci{#2}		
		\newcommand\drct{#3}		
		\foreach \k in {0,...,\reprn} {		
			\node at (-\drci / \reprp, \drcnn / \reprp * -\k / \drcnn) (p\drci\k) {\textbullet} ;
			\node[below right] at (p\drci\k) {$\repr{\drci}{\k}{\drct}$};
		} 
		\fill (-\drci / \reprp, -1.0 * \drcnn / \reprp) coordinate (p\drci\drcnn) rectangle +(1 / \reprp, \drcnn / \reprp * 0.2 / \drcnn); 	
		\draw (p\drci\drcnn) node[below right] (a\drci\drcnn) {$\repr{\drci}{\drcnn}{\drct}$} 
						rectangle +(1 / \reprp, 0.2 / \drcnn * \drcnn / \reprp); 	
	\end{scope}
} 

\newcommand{\drawreprxnn}[3][]{
	\begin{scope}[{#1}]
		\newcommand\drcnn{\the\numexpr\reprn+1\relax}	
		\newcommand\drci{#2}		
		\newcommand\drct{#3}		
		\fill (-\drci / \reprp, -1.0 * \drcnn / \reprp) coordinate (p\drci\drcnn) rectangle +(1 / \reprp, \drcnn / \reprp * 0.2 / \drcnn); 	
		\draw (p\drci\drcnn) node[below right] (a\drci\drcnn) {$\repr{\drci}{\drcnn}{\drct}$} 
						rectangle +(1 / \reprp, 0.2 / \drcnn * \drcnn / \reprp); 	
	\end{scope}
} 

\newcommand{\drawrepr}[2][]{
	\begin{scope}[{#1}]
		\newcommand\drt{#2}		
		\node at (0,0) (p00) {\textbullet} node[below right] (a00) {$\repr{0}{0}{\drt}$};
		\foreach \i in {1,...,\reprp}	{			
			\drawreprcolumn{\i}{\drt};
		} 
	\end{scope}
} 

\newcommand{\drawsamplex}[4][]{
	\begin{scope}[{#1}]
		\newcommand\dss{1 / \reprp}	
		\foreach \i in {{#2},...,{#3}} {
			\draw (-\i * \dss, 0) .. 
						controls (-\i * \dss + 0.25 * \dss, \dss + 1 * \dss * \dss - \i * \dss * \dss) 
						      and (-\i * \dss + 0.25 * \dss, \dss + 2 * \dss * \dss - \i * \dss * \dss) .. 
				(-\i * \dss+0.5 * \dss, 0) .. 
						controls +(0.25 * \dss, -\dss - 2 * \dss * \dss + \i * \dss * \dss) .. 
				(-\i * \dss + \dss, 0);
		} 
	\end{scope}
} 

\newtheorem{theorem}{Theorem}
\newtheorem{definition}{Definition}
\newtheorem{lemma}[theorem]{Lemma}
\newtheorem{rem}[theorem]{Remark}

\newcommand{\inte }{{\rm int}\,}

\newcommand{\cl }{{\rm cl}\,}

\newcommand{\comment}[1]{\mbox{}}

\def\qed{{\hfill{\vrule height5pt width5pt depth0pt}\medskip}}
\newcommand{\R}{\mathbb{R}}

\begin{document}

\begin{center}
{\bf \LARGE
Algorithm for rigorous integration of Delay Differential
Equations and the computer-assisted
proof of periodic orbits in the Mackey-Glass equation}

\vskip 0.5cm

{\large Robert Szczelina$^{1}$,  Piotr Zgliczy\'nski$^2$}

\vskip 0.5cm

$^1$ \textbf{Corresponding author}, Ma\l{}opolska~Center~of~Biotechnology, Jagiellonian~University,
Gronostajowa~7a, 30-374~Krak\'ow, Poland, e-mail:~robert.szczelina@uj.edu.pl, 

\vskip 0.25cm

$^2$ Institute~of~Computer~Science~and~Computational~Mathematics, Jagiellonian~University,
\L{}ojasiewicza~6, 30-348~Krak\'ow, Poland, e-mail:~zgliczyn@ii.uj.edu.pl

\vskip 0.5cm

\today

\vskip 0.5cm
\end{center}

\begin{abstract}
We present an algorithm for the rigorous integration of Delay Differential
Equations (DDEs) of the form $x'(t)=f(x(t-\tau),x(t))$. As an
application, we give a computer assisted proof of the existence of two
attracting periodic orbits (before and after the first period-doubling bifurcation)
in the Mackey-Glass equation.
\end{abstract}

\section{Introduction}
\label{sec:intro}


The goal of this paper is to present an algorithm for the rigorous  integration of Delay Differential Equations (DDEs) of
the form
\begin{equation}
\label{eq:dde}
\dot{x}(t) = f\left(x(t-\tau), x(t)\right), \quad x \in \R
\end{equation}
where $0 < \tau \in \mathbb{R}$.

Despite its apparent simplicity, Equation~\eqref{eq:dde} can generate
all kinds of possible dynamical behaviours: from simple stationary
solutions to chaotic attractors. For example, this happens for the well-known Mackey-Glass equation:
\begin{equation}
\label{eq:mackey-glass}
\dot{x}(t) = \beta \cdot {{x(t-\tau)} \over {1 + x^n(t-\tau)}} -\gamma \cdot x(t), \quad x \in \R,
\end{equation}
for which numerical experiments show the existence of a series of
period doubling bifurcations which lead to the creation of an apparent chaotic
attractor \cite{mackey-glass,mackey-glass-scholarpedia}. Later in the paper, we will apply our
rigorous integrator to this equation.

There are many important works that establish the existence and the
shape of a (global) attractor under various assumptions on $f$ in
Equation~\eqref{eq:dde}. Much is known about systems of the form
\mbox{$\dot{x} = -\mu x(t) + f\left(x(t-1)\right)$} when $f$ is strictly
monotonic, either positive or negative \cite{krisztin-review}. Let us mention here
a few developments in this direction.
\
Mallet-Paret and Sell used discrete Lyapunov
functionals to prove a Poincar\'e-Bendixson type theorem
for special kind of monotone systems \cite{mallet-paret-sell-bendixon}.
\
Krisztin, Walther and Wu have conducted a thorough study on systems
having a monotone positive feedback, including studies on the
conditions needed to obtain the shape of a global attractor, see
\cite{krisztin-book} and references therein.
\
In the case of a monotonic
positive feedback $f$ and under some assumptions on the stationary
solutions, Krisztin and Vas  proved that there exist large amplitude slowly oscillatory periodic
solutions (LSOPs) which revolve around more than one stationary
solution. Together with their unstable manifolds, connecting them with
the classical spindle-like structure, they constitute the full global
attractor for the system \cite{dde-lsop-krisztin}.
\
In a recent work, Vas showed that $f$ may be chosen such that the
structure of the global attractor may be arbitrarily complicated
(containing an arbitrary number of unstable LSOPs)
\cite{vas-big-attractor}.

\
Lani-Wayda and Walther were able to construct systems of the form
$\dot{x} = f\left(x(t-1)\right)$ for which they proved the existence of transversal homoclinic trajectory, and a hyperbolic set on which the dynamics are chaotic.
\cite{lani-wayda-walther-chaos}.

Srzednicki and Lani-Wayda  proved, by the use of the generalized Lefshetz fixed point
theorem, the existence of multiple periodic orbits and the existence of chaos
for some periodic, tooth-shaped (piecewise linear) $f$ \cite{srzednicki-wayda}.

The results from
\cite{dde-lsop-krisztin,srzednicki-wayda,lani-wayda-walther-chaos,vas-big-attractor},
while impressive, are established for functions which are close to
piecewise affine ones. The authors of these works
construct equations where an interesting behaviour appears,
however it is not clear how to apply their techniques for some well
known equations.

In recent years, there appeared many computer assisted proofs of various
dynamical properties for ordinary differential equations
and (dissipative) partial differential equations by an application of
 arguments from the geometric theory of dynamical systems
plus the rigorous integration, see for example
\cite{arioli-koch-ks,kapela-zgliczynski,cap-lorenz,tucker-lorenz,wilczak-map-1,
zgliczynski-pde-1}
and references therein. By the computer assisted proof we
understand a computer program which rigorously checks assumptions of
abstract theorems. This
paper is an attempt to extend this approach to the case of DDEs by
creating a rigorous forward-in-time integration scheme for
Equation~\eqref{eq:dde}. By the rigorous integration we
understand a computer procedure which produces rigorous
bounds for the true solution. In the case of DDEs, the integrator
should reflect the fact that, after the integration time longer than the
delay $\tau$, the solution becomes smoother, which gives the compactness of the
evolution operator. Having an integrator, one should be
able to directly apply standard tools from dynamics such as Poincar\'e
maps, various fixed point theorems, etc. In this paper, as an application, we present
computer-assisted proofs of the existence of two stable
periodic orbits for Mackey-Glass equation, however we do not prove
that these orbits are attracting.

There are several papers that deal with computer assisted proofs of
periodic solutions to DDEs \cite{lessard1, lessard2010, zalewski}, but
the approach used there is very different from our method. These works
transform the question of the existence of periodic orbits into a boundary
value problem (BVP), which is then
solved by using the Newton-Kantorovich theorem \cite{lessard1,lessard2010}
 or the local Brouwer degree  \cite{zalewski}. It is clear,
that the rigorous integration may be used to obtain more diverse
spectrum of results.  There are also several interesting results
that apply rigorous numerical computations
to solve problems for DDEs \cite{dde-wright-krisztin, dde-ricker-krisztin}, but they
do not rely on the rigorous, forward in time integration of DDEs.

The rest of the paper is organized as follows. Section~\ref{sec:integrator}
describes the theory and algorithms for the integration of
Equation~\eqref{eq:dde}. Section~\ref{sec:poinc-map} defines the
notion of the Poincar\'e map and discusses computation of the
Poincar\'e map using the rigorous integrator. Section~\ref{sec:results} presents an application of
the method to prove the existence of two stable periodic orbits in the Mackey-Glass
equation (Equation~\eqref{eq:mackey-glass}). Here, we investigate case for
$n = 6$ (before the first period doubling bifurcation) and for $n = 8$
(after the first period doubling bifurcation). To the best of our knowledge,
these are the first rigorous proofs of the existence of these orbits. Presented methods has been  also successfully used by the first author to prove 
the existence of multiple periodic
orbits in some other nonlinear DDEs \cite{szczelina-ejde}. 

\subsection{Notation}
\label{sub:notation}

We use the following notation. For a function $f :
\mathbb{R} \to \mathbb{R}$, by $f^{(k)}$ we denote the $k$-th derivative
of~$f$. By $f^{[k]}$ we denote the term $\frac{1}{k!}\cdot
f^{(k)}$. In the context of piecewise smooth maps by $f^{(k)}(t^-)$ and
$f^{(k)}(t^+)$ we denote the one-sided derivatives  $f$ w.r.t. $t$.

For $F : \mathbb{R}^m \to \mathbb{R}^n$ by $D F (z)$ we denote the matrix
$\left( \frac{\partial F_i}{\partial x_j}(z) \right)_{i \in \{1, .., n\}, j \in \{1, .., m\}}$.

For a given set $A$, by $\cl(A)$ and $\inte(A)$ we denote the closure and
interior of $A$, respectively (in a given topology e.g. defined by
the norm in the considered Banach space).

Let $A = \Pi_{i=1}^n [a_i, b_i]$ for $a_i \le b_i$, $a_i, b_i \in
\mathbb{R}$. Then, we call $A$ \emph{an interval set} (a product of closed
intervals in $\mathbb{R}^n$). For any $A \subset \mathbb{R}^n$
we denote by $hull(A)$ a minimal interval set, such that
$A \subset hull(A)$. If $A \subset \R$ is bounded then
$hull(A) = [\inf(A), \sup(A)]$.
For sets $A \subset \mathbb{R}$,
$B \subset \mathbb{R}$, $a \in \mathbb{R}$ and for some
binary operation $\diamond : \mathbb{R} \times \mathbb{R} \to \mathbb{R}$
we define $A \diamond B = \left\{a \diamond b: a \in A, b \in B\right\}$
and $a \diamond A = A \diamond a = \{ a \} \diamond A$.
Analogously, for $g : \mathbb{R} \to \mathbb{R}$ and a set
$A \in \mathbb{R}$ we define $g(A) = \{g(a) \ | \ a \in A\}$.

For $v \in \mathbb{R}^n$ by $\pi_i v$ for $i \in \{1, 2, .., n \}$ we
denote the projection of $v$ onto the $i$-th coordinate.
For vectors $u, v \in \mathbb{R}^n$ by $u \cdot v$ we
denote the standard scalar product:
$u \cdot v = \sum_{i=1}^{n} \pi_i v \cdot \pi_i u$

We denote by $C^r\left(D, \mathbb{R}\right)$ the space of all
functions of class $C^r$ over a compact set $D \subset \R$,
equipped with the supremum $C^r$ norm: $\| g \| = \sum_{i=0}^{r} \sup_{x \in D} | g^{(i)}(x) |$.
In case $D = [-\tau, 0]$, when $\tau$ is known from the context,
we will write $C^k$ instead of $C^k\left([-\tau,0], \R\right)$.

For a given function $x : [-1, a) \to \mathbb{R}$, $a \in \R_+ \cup \{ \infty \}$
for any $t \in [0, a)$ we denote by $x_t$ a function such that $x_t(s) = x(t+s)$
for all $s \in [-1, 0]$.

We will often use a symbol in square brackets, e.g.
$[r]$, to denote a set in $\R^m$.
Usually it will happen in formulas used
in algorithms, when we would like to
stress the fact that a given variable represents a set. If both
variables $r$ and $[r]$ are used simultaneously then
usually $r$ represents a value in $[r]$, however this is
not implied by default and it will be always
stated explicitly. Please note, that the notation $[r]$ does not
impose that the set $[r]$ is of any particular shape, e.g. an interval box.
We will always explicitly state if the set is an interval box.

For any set $X$ by $mid(X)$ we denote the midpoint of $hull(X)$
and by $diam(X)$ the diameter of $hull(X)$.

\subsection{Basic properties of solutions to DDEs}
For the convenience of the reader, we recall (without proofs) several classical results for DDEs \cite{dde-ode-book}.

We define the \emph{semiflow} $\varphi$ associated to Equation~\eqref{eq:dde} by:
\begin{equation}
\label{def:semiflow}
\varphi : \mathbb{R}_+ \times C^0([-\tau, 0], \R)   \ni (t,\psi) \mapsto x^\psi_t \in C^0([-\tau, 0], \R).
\end{equation}
where $x^\psi : [-\tau, a_\psi) \to \R$ is a solution to a Cauchy problem:
\begin{equation}
\label{eq:dde-ivp}
\begin{cases}
\dot{x} = f\left(x(t-\tau), x(t)\right), \quad t \ge 0, \\
x(t) = \psi(t), \quad t \in [-\tau, 0],
\end{cases}
\end{equation}
for a maximal $a_\psi \in \R_+ \cup \{ \infty \}$ such that
the solution exists for all $t < a_\psi$.

\begin{lemma}[Continuous (local) semiflow]
\label{lem:semiflow}
If $f$ is (locally) Lipshitz, then $\varphi$ is a (local) continuous semiflow on $C^0([-\tau, 0], \R)$.
\end{lemma}

\begin{lemma}[Smoothing property]
\label{lem:smoothing}
Assume $f$ is of class $C^m$, $m > 0$. Let $n \in \mathbb{N}$ be given and let $t \ge n \cdot \tau $.
If $x_0 \in C^0$ then $x_t = \varphi(t, x_0)$ is of class at least $C^{\min(m+1, n)}$.
\end{lemma}

The smoothing of solutions gives rise to some interesting objects
in DDEs \cite{walther-manifold}. Assume for a while that $f \in C^\infty$.
Then for any $n \geq 0$ there exists a set (in fact a manifold) $M^n \subset C^n$,
such that $M^n$ is forward invariant under $\varphi$.

It is easy to see that for $n=1$ we have:
\begin{equation*}
  M^1 := \left\{ x \in C^1 \ | \ x'(0^-)=f(x(-\tau),x(0)) \right\},
\end{equation*}
and the conditions for $M^n$ with $n > 1$ can be simply obtained by
differentiating both sides of \eqref{eq:dde}. We follow \cite{walther-manifold}
and we call $M^n$ a \emph{$C^n$ solution manifold}.

Notice that $M^n \subset M^k$ for $k \le n$ and
 $\varphi(k \tau, \cdot) : M^n \to M^{n+k}$.

\section{Rigorous  integration of DDEs}
\label{sec:integrator}
This section is a reorganized excerpt from the PhD dissertation of
the first author (Robert Szczelina). A detailed analysis of results
from numerical experiments with the proposed methods, more
elaborate description of the algorithms, and detailed pseudo-codes of the routines
can be found in the original dissertation \cite{szczelina-phd}.

\subsection{Finite representation of ,,sufficiently smooth'' functions}

Here, we would like to present the basic blocks used in the algorithm for the
rigorous integration of Equation~(\ref{eq:dde}). The idea is to
implement the Taylor method for Equation~\eqref{eq:dde} based on the piecewise
polynomial representation of the solutions plus a remainder term. We will work
on the equally-spaced grid and we will fix the step size of the Taylor method
to match the selected grid.

\begin{rem}
\label{rem:simplification}
In this section, for the sake of simplicity of
presentation, we assume that $\tau = 1$.
All computations can be easily redone for any delay $\tau$.

We also assume that r.h.s. $f$ of Equation~\eqref{eq:dde} is
,,sufficiently smooth'' for various expressions to make sense.
The class of $f$ in \eqref{eq:dde} restricts
the possible order of the Taylor method that can be used in our algorithms,
that is, if $f$ is  of class $C^n$, then we can use Taylor method of order at most $n$.
Therefore, thorough the paper it can be assumed that $f \in C^\infty$.
This is a reasonable assumption in the case of applications of computer-assisted proofs
where r.h.s. of equations are usually presented as a composition of elementary functions.
The Mackey-Glass equation \eqref{eq:mackey-glass} is a good example (away from $x = -1$).
\end{rem}

We fix two integers $n \ge 0$ and $p > 0$ and we set $h = \frac{1}{p}$.

\begin{definition}
By $C^{n}_{p}$ we denote the set of all functions
$g : [-1, 0] \to \mathbb{R}$ such that, for $i \in \left\{1,..,p\right\}$,
 we have:
\begin{itemize}
\item $g$ is $(n+1)$-times differentiable on $(-i \cdot h, -i \cdot h + h)$,

\item $g^{(k)}\left(-i \cdot h^+\right)$ exists for all $k \in \{0,..,n+1\}$ and $\lim\limits_{\xi \to 0^+} g^{(k)}(-i \cdot h+\xi) = g^{(k)}(-i \cdot h^+)$,

\item $g^{(n+1)}$ is continuous and bounded on $(-i \cdot h, -i \cdot h + h)$.
\end{itemize}
\end{definition}

From now on, we will abuse the notation and we will denote
the right derivative $g^{(k)}(-i \cdot h^+)$ by
$g^{(k)}(-i \cdot h)$ unless explicitly stated otherwise.
The same holds for $g^{[k]}(-i \cdot h^+)$.
Under this notation, it is clear that we can represent $g \in C^n_p$
by \emph{a piecewise Taylor} expansion on each
interval $\left[-i \cdot h, -i \cdot h + h\right)$.
For $t = - i \cdot h + \epsi$ and $1 \le i \le p$, $0 \le \epsi < h$
we can write:
\begin{equation}
\label{eq:representation-taylor-interpretation}
g(t) = \sum_{k=0}^{n} g^{[k]}(-i \cdot h) \cdot \epsi^{k} + g^{[n+1]}\left(-i \cdot h + \xi(\epsi)\right) \cdot \epsi^{n+1},
\end{equation}
with $\xi(\epsi) \in [0, h]$.

In our approach, we store the piecewise Taylor
expansion as a finite collection of coefficients
$g^{[k]}(- i \cdot h)$ and interval bounds on $g^{[n+1]}(\cdot)$ over the
whole interval $[-i \cdot h, -i \cdot h + h]$ for $i \in \{1,..,p\}$.
Our algorithm for the rigorous integration of \eqref{eq:dde} will then
produce rigorous bounds on the solutions to \eqref{eq:dde} for
initial functions defined by such piecewise Taylor expansion.

Please note, that we are using here a word \emph{functions} instead
of a single \emph{function}, as, because of the bounds on $g^{[n+1]}$ over intervals $[-i\cdot h, -i \cdot h + h]$,
the finite piecewise Taylor expansion describes an infinite set of
functions in general. This motivates the following definitions.

\begin{definition}
\label{def:representation}
Let $g \in C^n_p$ and let
$\mathbb{I} : \{ 1, .., p\} \times \{ 0, .., n\} \to \{1, .., p \cdot (n+1) \}$
be any bijection.

\emph{A minimal $(p,n)$-representation} of $g$ is a pair
$\bar{g} = (a, B)$ such that
\begin{eqnarray*}
a \in \R^{p \cdot (n+1) + 1} & \quad & B \subset \R^{p} \text{ is an interval set} \\
\pi_1 a = g(0) & \quad & \pi_i B = \left[\inf\limits_{\xi \in (0,h)} g^{[n+1]}\left(-i h + \xi\right), \sup\limits_{\xi \in (0,h)} g^{[n+1]}\left(-i h + \xi\right) \right]\\
\pi_{1+\mathbb{I}(i,k)} a = g^{[k]}(-i h^+) & \text{ for } & 1 \le i \le p,\ \ 0 \le k \le n
\end{eqnarray*}
\end{definition}

Please note, that the index function $\mathbb{I}$
 should be simply understood as an ordering  in which, during computations, we
store coefficients in a finite dimensional vector - its precise definition
is only important from the programming point of view, see Section~\ref{subsec:optimization} for a particular
example of $\mathbb{I}$.
So, in this
paper for theoretical considerations, we would like to use the following $\repr[g]{i}{k}{}$ notation
instead:
\begin{itemize}
\item $\repr[g]{0}{0}{} := \pi_{1} a$,

\item $ \repr[g]{i}{k}{} := \pi_{1+\mathbb{I}(i,k)} a$,

\item $\repr[g]{i}{n+1}{} := \pi_i B$
\end{itemize}
We call $\repr[g]{i}{k}{}$ \emph{the (i,k)-th coefficient of the
representation} and $\repr[g]{i}{n+1}{}$ \emph{the i-th remainder
of the representation}. The interval set $B$ is called
\emph{the remainder of the representation}. We will call the
constant $m = p \cdot (n+2) + 1$ the \emph{size of the $(p,n)$-representation}.
When parameters $n$ and $p$ are known from the context we will omit
them and we will call $\bar{g}$ \emph{the minimal representation}
of $g$.

\begin{definition}
\label{def:support}
We say that $G \subset C^n_p$ is a \emph{$(p,n)$-f-set} (or $(p,n)$-functions set)
if there exists bounded set $[\bar{g}] = (A, C) \subset \R^{p\cdot(n+1)+1} \times \R^p  = \R^m$
such that
\begin{equation*}
G = \left\{f \in C^{n}_{p} \ | \ \bar{f} \subset [\bar{g}] \mbox{ for the minimal $(p,n)$-representation } \bar{f} \mbox{ of $f$} \right\}.
\end{equation*}
As the set $[\bar{g}]$ contains the minimal representation of $f$ for any
$f \in G$, we will also say that $[\bar{g}]$ is a \mbox{$(p,n)$-representation}
of $G$. We will also use $G$ and $[\bar{g}]$ interchangeably and we will
write $f \in [\bar{g}]$ for short, if the context is clear.
\end{definition}
Please note that the minimal $(p,n)$-representation $\bar{g}$ of $g$
defines  $(p,n)$-f-set $G \subset C^n_p$, which, in general, contains
more than the sole function $g$. Also, in general, for any $(p,n)$-f-set $G$
there are  functions $g \in G$ which are discontinuous at
grid points $-i \cdot h$ (see \eqref{eq:representation-taylor-interpretation}).
Sometimes we will need to assume higher regularity, therefore we define:
\begin{definition}
\label{def:c-k-support}
Let $G=[\bar{g}]$ be a $(p,n)$-f-set.
\emph{The $C^k$-support of $G$} is defined as:
\begin{equation*}
Supp^{(k)}([\bar{g}]) \ := \ Supp^{(k)}(G) \ := \ G \cap C^k.
\end{equation*}
For convenience we also set:
\begin{equation*}
Supp([\bar{g}]) \ := \ Supp(G) \ := \ G.
\end{equation*}
\end{definition}
Please note that $C^n_p \supset Supp(G) \ne Supp^{(0)}(G) \subset C^0$. It may also happen that $Supp^{(k)}(G)=\emptyset$ for nonempty $G$ even for $k=0$.

Now we present three simple facts about the convexity of the support sets. These
properties will be important in the context of the computer assisted proofs
and in an application of Theorem~\ref{thm:schauder} to $(p,n)$-f-sets in Section~\ref{sec:results}.
\begin{lemma}
For $(p,n)$-representations $[\bar{g}], [\bar{f}] \subset \mathbb{R}^m$,
$m = p \cdot (n+2) + 1$, the following statements hold true:
\begin{itemize}
\item If $[\bar{g}] \subset [\bar{f}]$, then
$Supp([\bar{g}]) \subset Supp([\bar{f}])$.

\item If $[\bar{g}]$ is a convex set in $\R^m$,
then $Supp([\bar{g}])$ is a convex set in $C^n_p$.

\item If $[\bar{g}]$ is a convex set in $\R^m$, then
$Supp([\bar{g}]) \cap C^k$ is a convex set for any $k \ge 0$.
\end{itemize}
\end{lemma}
We omit the easy proof.

To extract information on
$g^{[k]}\left(-{i \over p} + \epsi\right)$ for any $i$ and $k$ having
only information stored in a $(p,n)$-representation, we introduce the following definition.

\begin{definition}
\label{def:cksymb}
Let $(p,n)$-representation $\bar{g}$ be given. We define
\begin{equation*}
c^{i,[k]}_{\bar{g}}(\epsi)  =  \sum_{l=k}^{n+1} {l \choose k} \cdot \epsi^{l-k} \cdot \repr[g]{i}{l}{},
\end{equation*}
for $0 < \epsi < \frac{1}{p}$, $1 \le i \le p$ and $0 \le k \le n+1$.
\end{definition}
We will omit subscript $\bar{g}$ in $c^{i,[k]}_{\bar{g}}(\epsi)$ if
it is clear from the context. The following lemma follows immediately
from the Taylor formula, so we skip the proof:
\begin{lemma}
\label{lem:c-k}
Assume $g \in C^n_p$ and its $(p,n)$-representation
$\bar{g}$ are given.
Then
for $0 < \epsi < \frac{1}{p}$, $1 \le i \le p$ and $0 \le k \le n+1$
\begin{equation*}
g^{[k]}\left(-\frac{i}{p}+\epsi\right) \in c^{i,[k]}\left(\epsi\right)
\end{equation*}
holds.
\end{lemma}

\comment{
\textbf{Proof:}
Let first look at the Taylor expansion of
$g\left(-\frac{i}{p} + \epsi \right)$:
\begin{equation*}
g\left(-\frac{i}{p} + \epsi \right)  =  g \left( -\frac{i}{p} \right) + \sum_{l=1}^{n} \frac{\epsi^l}{l!} \cdot g^{(l)} \left( -\frac{i}{p} \right) + \frac{\epsi^{n+1}}{(n+1)!} \cdot g^{(n+1)} \left( -\frac{i}{p}+\xi_0 \right),
\end{equation*}
where $\xi_0 = \xi_0(\epsi) \in \left[0, \epsi \right]$. Using
Definition~\ref{def:representation} of the representation we
obtain that
\begin{equation*}
g^{[0]}\left(-\frac{i}{p} + \epsi\right) = g\left(-\frac{i}{p} + \epsi\right) \ \in \ \sum_{l=0}^{n+1} \epsi^l \cdot \repr[g]{i}{l}{} \ = \ c^{i,[0]}(\epsi).
\end{equation*}
Now, similarly, we can express each
$g^{(k)}\left(-\frac{i}{p} + \epsi\right)$ by:
\begin{eqnarray*}
g^{(k)}\left(-\frac{i}{p} + \epsi \right) & = & \sum_{l=k}^{n} \epsi^{l-k} \cdot \frac{g^{(l)} \left( -\frac{i}{p} \right)}{(l-k)!} + \\
 & + & \epsi^{n+1-k} \cdot \frac{g^{(n+1)} \left( -\frac{i}{p}+\xi_k \right)}{(n+1-k)!},
\end{eqnarray*}
for some $\xi_k = \xi_k(\epsi) \in \left[0, \epsi \right]$. Dividing
both sides by $k!$, setting $g^{(l)} = l! \cdot g^{[l]}$ and
estimating $\xi_k$ by the representation we get the following:
\begin{equation*}
g^{[k]}\left(-\frac{i}{p}+\epsi\right) \ \in \ \sum_{l=k}^{n+1} \frac{l!}{k! \cdot (l-k)!} \cdot \epsi^{l-k} \cdot \repr[g]{i}{l}{} \ = \ c^{i,[k]}(\epsi).
\end{equation*}
\qed
}

Before proceeding to the presentation of the integration procedure, we would
like to discuss the problem of obtaining Taylor coefficients
of a solution $x$ to Equation~\eqref{eq:dde} at a given time $t$
(whenever they exist).
From Equation~\eqref{eq:dde}, we have
(we remind that, at grid points, by the derivative we
mean the right derivative):
\begin{eqnarray*}
x^{(k)}(t) & = & \frac{d^{k-1}}{dt^{k-1}}f(x(t-1),x(t)).
\end{eqnarray*}
For example, in case of $k=1$, we obviously have:
\begin{equation*}
x^{(1)}(t) = f\left(x(t-1), x(t)\right),
\end{equation*}
and in case $k=2$, by applying the chain rule, we get:
\begin{eqnarray*}
x^{(2)}(t) & = & \frac{\partial f}{\partial z_1} \left(x(t-1), x(t)\right) \cdot x^{(1)}(-1) +  \\
& + & \frac{\partial f}{\partial z_2} \left(x(t-1), x(t)\right) \cdot x^{(1)}(0)
\end{eqnarray*}
If we define a function $F_{(1)} : \R^4 \to \R$ as
\begin{eqnarray*}
F_{(1)}(z_1, z_2, z_3, z_4) & = & \frac{\partial f}{\partial z_1} \left(z_1, z_3\right) \cdot z_2 + \frac{\partial f}{\partial z_3} \left(z_1, z_3\right) \cdot z_4,
\end{eqnarray*}
then we see that
\begin{eqnarray*}
x^{(2)}(t) & = & F_{(1)}\left(x(t-1), x^{(1)}(t-1), x(t), x^{(1)}(t)\right).
\end{eqnarray*}
Now, by a recursive application of the chain rule, we can obtain a  family of functions
$\mathbb{F}_f = \left\{F_{(k)} : \mathbb{R}^{2 \cdot (k+1)} \to \mathbb{R}\right\}_{k \in \mathbb{N}}$
such that:
\begin{eqnarray*}
x^{(k+1)}(t) = F_{(k)}\left(x(-1+t), .., x^{(k)}(-1+t), x(t), .., x^{(k)}(t)\right).
\end{eqnarray*}
By setting
\begin{equation}
\label{eq:F[k]}
F_{[k]}(z_1, .., z_{2\cdot (k+1)}) = \frac{1}{k!} F_{(k)}\left(0! \cdot z_1, .., k! \cdot z_{k+1}, 0! \cdot z_{k+2}, .., k! \cdot z_{2 \cdot (k+1)}\right),
\end{equation}
we can write similar identity in terms of the Taylor
coefficients $x^{[k]}$:
\begin{equation}
\label{eq:F[k]=x[k]}
x^{[k+1]}(t) = \frac{1}{k+1} \cdot F_{[k]}\left(x^{[0]}(-1+t), .., x^{[k]}(-1+t), x^{[0]}(t), .., x^{[k]}(t)\right).
\end{equation}
As we are using the Taylor coefficients instead of derivatives to represent
our $(p,n)$-f-sets, this notation would be more suitable to describe computer algorithms.
From now on, we will also slightly abuse the notation and we will denote $F_{[k]}$
by $F^{[k]}$ and $F_{(k)}$ by $F^{(k)}$. This is reasonable, since, for
a function $F : \R \to \R$ defined by $F(t) := f(x(t-1), x(t))$, we have:
\begin{eqnarray*}
F^{(k)}(t) &  =  & F_{(k)}\left(x^{(0)}(-1+t), .., x^{(k)}(-1+t), x^{(0)}(t), .., x^{(k)]}(t)\right), \\
           & and & \\
F^{[k]}(t) &  =  & F_{[k]}\left(x^{[0]}(-1+t), .., x^{[k]}(-1+t), x^{[0]}(t), .., x^{[k]}(t)\right).
\end{eqnarray*}

\begin{rem}
The task of obtaining family $\mathbb{F}_f$ by directly and analytically
applying the chain rule may seem quite tedious, especially, if one will be
required to supply this family as implementations of computer procedures.
It turns out, that this is not the case for a wide class of functions.
In fact, only the r.h.s. of Equation~\eqref{eq:dde} needs to be implemented
and the derivatives may be obtained by the means of the  automatic
differentiation (AD) \cite{AD1,AD2}. We use Taylor coefficients
$x^{[k]}$ to follow the notation and implementation of AD in the CAPD
library \cite{www-capd} which provide a set of rigorous interval
arithmetic routines used in our programs.
\end{rem}

\subsection{One step of the integration with fixed-size step $h = \frac{1}{p}$}

We are given $(p,n)$-f-set $\bar{x}_0$ and the task is to obtain
$\bar{x}_h$ - a \mbox{$(p,n)$-f-set}
such that $\varphi(h, \bar{x}_0) \subset \bar{x}_h$.
We will denote the procedure of computing $\bar{x}_h$ by $I_h$, that is:
\begin{equation*}
\bar{x}_h = I_h(\bar{x}_0).
\end{equation*}

First of all, we consider how $\bar{x}_0$ and $\bar{x}_h = I_h(\bar{x}_0)$
relate to each other. Their mutual alignment  is shown in
Figure~\ref{fig:integrator}.

\begin{figure}
    \newcommand\reprp{4}
    \newcommand\reprn{2}
    \newcommand\reprnn{3}
    \centering
    \begin{tikzpicture}[scale=1.4*\reprp]
		\coordinate (upshift) at (1 / \reprp,  1 + 0.5 / \reprn);		
		\coordinate (downshift) at (0, -0.65 / \reprn);		
		\draw[->] (-1-1.0 / \reprp, 0) -- (2.5 / \reprp, 0) node[right] (ttt) {$t$}; 
		\drawsamplex[dotted]{-1}{\the\numexpr\reprp+1\relax};

		\begin{scope}
			\draw[color=lightgray] (0, 0) -- (0, 0);
			\foreach \i in {0,...,\reprp} { 
				\reprnode[shift={(downshift)}]{\i}{\reprnn}{p}
				\draw[color=lightgray] (-\i / \reprp, 0) -- (p); 
			} 

			\drawrepr[active, shift={(downshift)}]{0};
		\end{scope}

		\begin{scope}[shift={(1.0 / \reprp,0)}]
			\draw[color=lightgray] (0, 0) -- (0, 0);
			\foreach \i in {0,...,\reprp} { 
				\reprnode[shift={(upshift)}]{\i}{0}{p}
				\draw[color=lightgray] (-\i / \reprp, 0) -- (p); 
			} 

			\drawrepr[active, shift={(upshift)}]{h};
		\end{scope}

		\foreach \k in {1,...,\reprn} {
			\reprnode[shift={(upshift)}]{1}{\k}{p}
			\fill[color=white] (p) circle (0.02);
			\draw (p) circle (0.02);
		}
		\reprnode[shift={(upshift)}]{0}{0}{p}
		\fill[color=white] (p) circle (0.03);
		\draw (p) circle (0.03);
		\draw (p) circle (0.02);
		\reprnode[shift={(upshift)}]{1}{\reprnn}{p};
		\draw[fill=white] (p) rectangle +(1 / \reprp, 0.2 / \reprnn); 
		\draw[fill=white] (p)+(0.1 / \reprp / \reprp, 0.1 / \reprp / \reprp) rectangle +(1 / \reprp - 0.1 / \reprp / \reprp, 0.2 / \reprnn - 0.1 / \reprp / \reprp); 

	\end{tikzpicture}
    \caption{\label{fig:integrator}
        A graphical presentation of the integrator scheme.
        We set $n=\reprn$ and $p=\reprp$.
        A $(p,n)$-representation is depicted as dots at grid points and rectangles
        stretching on the whole intervals between consecutive grid points.
        The dot is used to stress the fact that
        the corresponding coefficient represents \emph{the value at a given grid point}.
        Rectangles are used to stress the fact that remainders are
        \emph{bounds for derivative over whole intervals}.
        Below the time line we have an initial $(p,n)$-representation.
        Above the time line we see a representation after
        one step of size $h = \frac{1}{p}$. Black solid
        dots and grey rectangles represent the values we do not need to compute -
        this is \emph{the shift part}. \emph{The forward part}, i.e. the elements
        to be computed, are presented as empty dots and an empty rectangle.
        The doubly bordered dot represents the
        exact value of the solution at the time $t = h = \frac{1}{p}$
        (in practical rigorous computations it is an interval bound on the value).
	   The doubly bordered empty rectangle is an enclosure for the $n+1$-st
	   derivative on the interval $\left[0, h\right]$.
    }
\end{figure}

We see that $\repr{i}{k}{h}$ overlap with
$\repr{i-1}{k}{0}$, so they can be simply \emph{shifted} to the new
representation - we call this procedure the \emph{Shift Part}.
Other coefficients need to be estimated using the dynamics generated
by Equation~\eqref{eq:dde}. We call this procedure
the \emph{Forward Part}. This procedure will be divided into three
subroutines:
\begin{enumerate}
	\item computing coefficients $\repr{1}{k}{h}$ for
           $k \in \{1, .., n\}$,
           \label{int-procedure-1}

	\item computing the remainder $\repr{1}{n+1}{h}$,
           \label{int-procedure-2}

	\item computing the estimate for $x_h(0)$
           (stored in $\repr{0}{0}{h}$).
           \label{int-procedure-3}
\end{enumerate}

\vskip 1em \noindent \textbf{Forward Part - Subroutine~\ref{int-procedure-1}}

This procedure is immediately obtained
by a recursive application of Equation~\eqref{eq:F[k]=x[k]}:
\begin{eqnarray*}
\repr{1}{0}{h} & = & \repr{0}{0}{0} \\
\repr{1}{k}{h} & = & \frac{1}{k} \cdot F^{[k-1]}\left(\repr{p}{0}{0}, .., \repr{p}{k-1}{0}, \repr{1}{0}{h}, .., \repr{1}{k-1}{h}\right),
\end{eqnarray*}
where $1 \le k \le n$.

\vskip 1em \noindent \textbf{Forward Part - Subroutine~\ref{int-procedure-2}}

This subroutine can be derived from
the Mean Value Theorem. We have for $\epsi < h$:
\begin{eqnarray*}
\frac{1}{(n+1)!} \cdot x^{(n+1)}(\epsi) & = & \frac{1}{(n+1)!} \cdot x^{(n+1)}(0) + \frac{1}{(n+1)!} \cdot x^{(n+2)}(\xi) \cdot \epsi =  \notag \\
\label{eq:compute-remainder} & = & \frac{1}{(n+1)} \cdot F^{[n]}\left(x(-1), .., x^{[n]}(-1), x(0), .., x^{[n]}(0), \right) \ + \label{eq:remainder} \\
& + & F^{[n+1]}\left(x(-1+\xi), .., x^{[n+1]}(-1+\xi), x(\xi), .., x^{[n+1]}(\xi)\right) \cdot \epsi \notag
\end{eqnarray*}
for some $0 \le \xi \le \epsi$. Let us look at the two terms
that appear on the r.h.s. of Equation~\eqref{eq:compute-remainder}.
The question is: can we estimate them by having only  $\bar{x}_0$ and already
computed $\repr{1}{k}{h}$ from Subroutine~\ref{int-procedure-1}?
Let us discuss each of these terms separately.

By Lemma~\ref{lem:c-k} we have for $0 \le k \le n+1$:
\begin{eqnarray*}
x^{[k]}(-1+\xi) & \in & c^{p,[k]}_{\bar{x}_0}\left([0, h]\right)
\end{eqnarray*}
Moreover by Definition~\ref{def:representation}  we know that:
\begin{eqnarray*}
x^{[k]}(-1) & \in & \repr{p}{k}{0} \\
x^{[k]}(0) & \in & \repr{1}{k}{h}
\end{eqnarray*}
So those terms can be easily obtained.
The problem appears when it comes to $x(\xi)$ and
$x^{[k]}(\xi)$, for $\xi \in (0,h)$

Assume for a moment that we have some a priori estimates
for $x([0,h])$ i.e. a set $Z \subset \R$ such that $x(\left[0, h]\right) \subset Z$.
We call this set the \emph{rough enclosure} of $x$ on the interval $[0, h]$.
Having  rough enclosure $Z$, we could apply Equation~\eqref{eq:F[k]=x[k]} (as
in the case of Subroutine~\ref{int-procedure-1}) to obtain the estimates on
$x^{[k]}\left([0, h]\right)$ for $k > 0$. So the question is: how to find a candidate $Z$
and prove that $x(\left[0, h]\right) \subset Z$?
The following lemma gives a procedure to test the later.
\begin{lemma}
\label{lem:rough-enclosure}
Let $Y \subset \mathbb{R}$ be a closed interval and let $x_0$ be a function defined
on $[-1, 0]$. Assume that the following holds true:
\begin{equation}
\label{eq:rough-enclosure}
Z := x_0(0) + \left[0, h\right] \cdot f\left(x_0\left(\left[-1,-1+h\right]\right), Y\right) \subset int(Y).
\end{equation}
Then the solution $x(t)$ of Equation~\eqref{eq:dde} with the initial condition
$x_0$ exists on the interval $\left[0, h\right]$ and
\begin{equation*}
x\left(\left[0,h\right]\right) \subset Z.
\end{equation*}
\end{lemma}
\textbf{Proof: } We can treat equation~\eqref{eq:dde} on the
interval $\left[0, h\right]$ as a non-autonomous ODE of the form:
\begin{equation*}
x' = f(a(t), x),
\end{equation*}
where $a(t) = x(t)$ for $t \in [-1, -1+h]$ is a known function. Now
the conclusion follows from the proof of the analogous theorem for
ODEs. The proof can be found in \cite{capd-zgliczynski}.
\qed

Using Lemma~\ref{lem:rough-enclosure}, a heuristic iterative algorithm may
be designed such that it starts by guessing an initial $Y$ and then
it applies Equation~\eqref{eq:rough-enclosure} to obtain $Z$.
In a case of failure of the inclusion, i.e. $Z \not \subset Y$, a bigger
$Y$ is taken in the next iteration. Please note, that this iteration may never stop
or produce unacceptably big $Y$, especially when the step-size $h$ is large.
Finding a rough enclosure is the only place in the  algorithm of the
integrator that can in fact fail to produce any estimates.
In such a case we are not able to proceed with the integration and we
signalize an error.

Now we can summarize the algorithm for Subroutine~\ref{int-procedure-2} as follows:
\begin{eqnarray*}
c^{[k]} & := & c^{p,[k]}([0, h]), \quad k \in \{0,..,n+1\} \\
d^{[0]} & := & Z \textrm{ as in equation \eqref{eq:rough-enclosure}}, \\
d^{[k]} & := & \frac{1}{k} \cdot F^{[k-1]}\left(c^{[0]}, .., c^{[k-1]}, d^{[0]}, .., d^{[k-1]}\right), \quad k \in \{1,..,n+1\} \\
a^* & := & \frac{1}{(n+1)} \cdot F^{[n]}\left(\repr{p}{0}{0}, .., \repr{p}{n}{0}, \repr{1}{0}{h}, .., \repr{1}{n}{h}\right) \\
b^* & := & F^{[n+1]}(c^{[0]}, .., c^{[n+1]}, d^{[0]}, .., d^{[n+1]}), \\
\repr{1}{n+1}{h} & := & a^* + b^* \cdot [0, h]
\end{eqnarray*}


\begin{rem}
Please note that the term $a^*$ is computed the same way as other
coefficients in Subroutine~\ref{int-procedure-1} and the rough enclosure
do not influence this term. In fact this is the $n+1$-st derivative of
the flow w.r.t. time. It is possible to keep track of those coefficients
during the integration and after $p$ steps (full delay) those coefficients may be used to
build a $(p,n+1)$-representation of the solutions - this is a direct
reflection of consequences of Lemma~\ref{lem:smoothing}.

This fact is also important for the compactness of the evolution
operator - an essential property that allows for an application of the topological fixed
point theorems in infinite dimensional spaces.
\end{rem}

\vskip 1.0em \noindent \textbf{Forward Part - Subroutine~\ref{int-procedure-3}}

The last subroutine of the forward part can be simply obtained by
using Definition~\ref{def:representation}
and Equation~\eqref{eq:representation-taylor-interpretation}:
\begin{equation*}
\repr{0}{0}{h} := \sum_{k=0}^{n}\repr{1}{k}{h} \cdot h^k + \repr{1}{n+1}{h} \cdot h^{n+1}.
\end{equation*}
Notice, that the possible influence of the usually over-estimated
rough enclosure $Z$ is present only in the
last term of the order $h^{n+1}$ so, for small $h$ (large enough $p$), it should not
be a problem.

\vskip 1em \noindent \textbf{The integrator - altogether}

Strictly speaking, the mapping $I_h$ does not produce a $(p,n)$-f-set
which exactly represents $\bar{x}_h = \varphi(h, \bar{x}_0)$. Instead, it returns
some bigger set $[\bar{x}_h]$ such that $\bar{x}_h$ is contained
in it. Of course, we are interested in obtaining a result as
close as possible to the set of true solutions represented
by $\varphi(h, \bar{x}_0)$. So, for technical reasons which
will be apparent in Section~\ref{subsec:wrapping}, we
decompose $I_h$ into
$I_h = \Phi + R$ such that $\Phi : \R^m \to \R^m$ and
$R : \mathcal{P}(\R^m) \to \mathcal{P}(\R^m)$.
Let $\Phi(\bar{x}) = \bar{\phi}$ and put:
\begin{eqnarray}
\repr[\phi]{i}{k}{} & := & \repr{i-1}{k}{}, \quad i \in \{2,..,p\}, k \in \{0,..,n\}, \label{eq:phi-start} \\
\repr[\phi]{1}{0}{} & := & \repr{0}{0}{}, \\
\repr[\phi]{1}{k}{} & := & \frac{1}{k} \cdot F^{[k-1]}\left(\repr{p}{0}{}, .., \repr{p}{k-1}{}, \repr[\phi]{1}{0}{}, .., \repr[\phi]{1}{k-1}{}\right), \quad k \in \{1,..,n\}, \\
\repr[\phi]{0}{0}{} & := & \sum_{k=0}^{n}\repr[\phi]{1}{k}{} \cdot h^k. \label{eq:phi-end}
\end{eqnarray}
Let $R(\bar{x}) = \bar{r}$ and put:
\begin{eqnarray*}
\repr[r]{i}{n+1}{} & := & \repr{i-1}{n+1}{}, \quad i \in \{2,..,p\}, \\
\repr[r]{1}{n+1}{} & := & a^* + b^* \cdot [0, h], \\
\repr[r]{0}{0}{} & := & h^{n+1} \cdot \left(a^\star + b^* \cdot [0, h]\right).
\end{eqnarray*}
The map $\Phi$ is called \emph{the Taylor part}, while the map $R$
is called \emph{the Remainder part}. This decomposition is important for an efficient
reduction of negative effects caused by using interval arithmetic,
primarily the wrapping effect, but also the
dependency problem to some extent.

\subsection{Reducing the wrapping effect}
\label{subsec:wrapping}

A representation of objects in computations  as the  interval sets has
its drawbacks. Possibly the most severe of them are the phenomena
called \emph{the wrapping effect} and \emph{the dependency problem}. Their influence is
so dominant that they are discussed in virtually every
paper in the field of rigorous computations
(see \cite{capd-zgliczynski} and references therein).
The dependency problem arises in interval arithmetic when two values
theoretically representing the same (or dependent) value are combined. The
most trivial example is an operation $x-x$ which is always $0$, but it is not the case for
intervals. For example, applying the operation to the interval $x = [-1,1]$
gives as the result the interval $[-2,2]$ which contains $0$ but it is far
bigger than we would like it to be.

\begin{figure}
    \centering{
        \includegraphics[width=80mm]{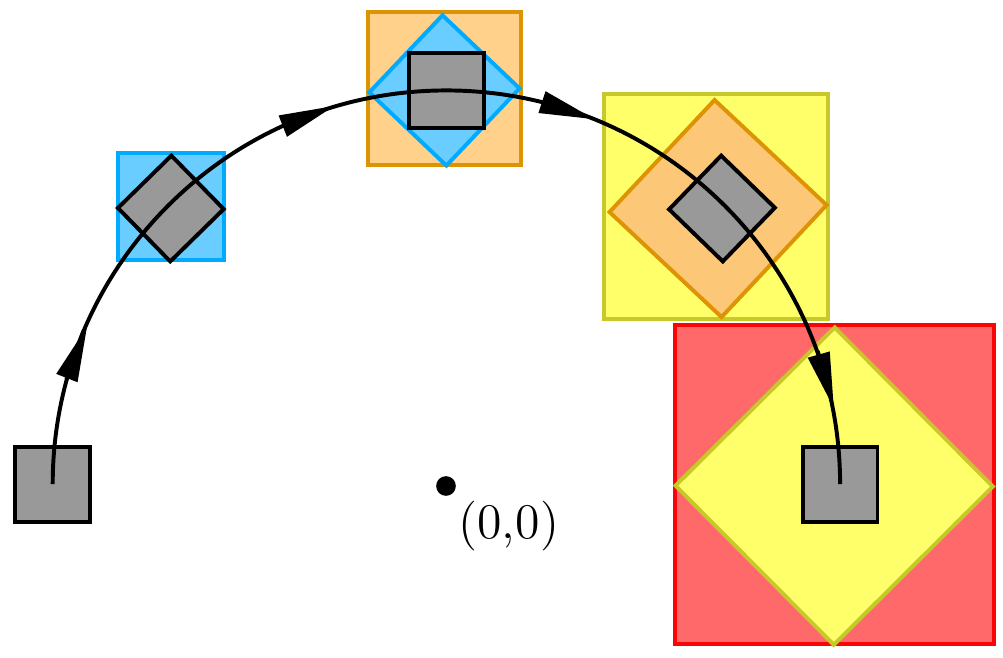}
    }

    \caption{\label{fig:wrapping-effect}
	An illustration of the wrapping effect problem for a classical,
	idealized mathematical pendulum ODE $\ddot{x} = -x$. The picture
	shows a set of solutions in a phase space $(x, \dot{x})$. The grey
     boxes present points of initial box
	moved by the flow. The colored boxes present the wrapping effect
	occurring at each step when we want to enclose the moving points in a
	product of intervals in the basic coordinate system. For example, the
     blue square on the left encloses the image of the first iteration. Its image is
     then presented with blue rhombus which is enclosed again by an orange
     square. Then the process goes on.
    }
\end{figure}

The wrapping effect arises when one intends to represent a result of some evaluation on sets
as a simple interval set. Figure~\ref{fig:wrapping-effect} illustrate this when we consider the rotation of the square.

One of the mostly used and efficient methods for reducing the impact of the wrapping effect and
the dependency problem was proposed by Lohner \cite{capd-lohner}. In the context of the iteration of
maps and the integration of ODEs, he proposed to represent sets by parallelograms,
i.e. interval sets in other coordinate systems. In the sequel we follow \cite{capd-zgliczynski}
and we sketch the Lohner methods briefly.

By $J$ we denote a computation of $I_h(\cdot)$ using point-wise evaluation of the Taylor part, i.e.:
\begin{equation*}
J([x]) := \left(\bigcup\limits_{x \in [x]} \Phi(x) \right) + R([x]).
\end{equation*}
Let us consider an iteration:
\begin{equation*}
[x_k] = J\left([x_{k-1}]\right), \quad k \in \mathbb{N}
\end{equation*}
with initial set $[x_0]$.

Let us denote $x_k = mid([x_k])$ and $[r_k] = [x_k] - x_k$.
By a simple argument based on the Mean Value Theorem \cite{capd-zgliczynski}, it can be shown that:
\begin{equation*}
[x_{k+1}] \subset \Phi(x_k) + \left[D\Phi\left([x_k]\right)\right] \cdot [r_k] + R([x_k]).
\end{equation*}
We can reformulate the problem of computing $[x_{k+1}]$ to the following
system of equations:
\begin{eqnarray}
[A_k] &=& \left[D\Phi\left([X_k]\right)\right], \label{eq:defAk} \\
x_{k+1} & = & mid\left(   {[\Phi(x_k) + R([X_k])]}   \right),  \\
{[z_{k+1}]} & = & \Phi({[x_k]}) + R([x_k]) - x_{k+1},  \\
{[r_{k+1}]} & = & [A_k] \cdot [r_k] + [z_{k+1}]. \label{eq:lohner-part}
\end{eqnarray}

Now the
reduction of the wrapping effect could be obtained by choosing
suitable representations of sets $[r_k]$ and a careful evaluation of
Equation~\eqref{eq:lohner-part}. The terminology used for this in
\cite{capd-zgliczynski} is \emph{the rearrangement computations}. We will
briefly discuss possible methods of handling Eq.~\eqref{eq:lohner-part}.

\noindent\textbf{Method 0} (Interval Set): Representation of $[r_k]$ by
an interval box and the direct evaluation of \eqref{eq:lohner-part} is
equivalent to directly computing $I_h([x_{k+1}])$. This method is called
\emph{an interval set} and is the least effective.

\noindent\textbf{Method 1} (Parallelepiped): we require that
$[r_k] = B_k \cdot [\tilde{r}_k]$ for $[\tilde{r}_k]$ being an
interval box and $B_k$ being an invertible matrix. Then
\eqref{eq:lohner-part} becomes:
\begin{equation*}
{[r_{k+1}]} = B_{k+1} \left(B_{k+1}^{-1}[A_k]B_{k+1} \cdot [r_k] + B_{k+1}^{-1}[z_{k+1}]\right).
\end{equation*}
Since it is difficult to obtain the exact matrix inverse in
computer calculations we will use interval matrices $[B_k]$ and $[B_k^{-1}]$
that contain $B_k$ and $B_k^{-1}$, respectively. Thus the equation on $[\tilde{r}]$ becomes:
\begin{equation*}
{[\tilde{r}_{k+1}]} = \left([B_{k+1}^{-1}][A_k][B_{k+1}]\right) \cdot [\tilde{r}_k] + [B_{k+1}]^{-1}[z_{k+1}].
\end{equation*}
If $[B_k]$'s are well-chosen, then the formula in brackets can be
evaluated to produce a matrix close to identity with very small diameter, thus the wrapping
effect reduction is achieved. The Parallelepiped method is obtained
when $B_{k+1}$ is chosen such that $B_{k+1} \in [A_k][B_{k+1}]$. This
approach is of limited use because of the need to compute the matrix inverse of a
general matrix $B_{k+1}$, which may fail or produce unacceptable
results if $B_{k+1}$ is ill-conditioned.

\noindent\textbf{Method 2} (Cuboid): this is a modification of
Method~1. In this method, we choose $U \in [A_k][B_{k+1}]$ and we do
the floating point approximate QR decomposition of $U = Q \cdot R$, where
$Q$ is close to an orthogonal matrix. Next we obtain matrix $[Q]$ by
applying the interval (rigorous) Gram-Schmidt method to $Q$, so there
exist orthogonal matrix $\tilde{Q} \in [Q]$ and
$Q^{-1} = Q^{T} \in [Q]^T$. We set $B_{k+1} = [Q]$,
$B_{k+1}^{-1} = [Q]^{T}$.

\noindent\textbf{Method 3} (Doubleton): this representation is
used in our computations as it proved to be the most efficient in
numerical tests \cite{szczelina-phd} and in other applications, see
\cite{capd-zgliczynski} and references therein. The original
idea by Lohner is to separate
the errors introduced due to the large size of initial data and
the local errors introduced by the numerical method at every
step. Namely we set:
\begin{eqnarray*}
{[r_{k+1}]} & = & {[E_{k+1}]}{[r_0]} + {[\tilde{r}_{k+1}]} \\
{[\tilde{r}_{k+1}]} & = & {[A_k]}{[\tilde{r}_k]} + {[z_{k+1}]} \\
{[E_{k+1}]} = {[A_k]}{[E_k]} & \quad & E_0 = Id
\end{eqnarray*}
where $[\tilde{r}_k]$ is evaluated by any method 0-2.
To reduce the possible wrapping effect in the product ${[A_k]}{[E_k]}$,
Lohner proposed the following:
\begin{eqnarray*}
{[r_{k+1}]} & = & C_{k+1}{[r_0]} + {[\tilde{r}_{k+1}]} \\
{[\tilde{r}_{k+1}]} & = & {[A_k]}{[\tilde{r}_k]} + {[z_{k+1}]} + \left({[A_k]}{C_k} - C_{k+1}\right) {[r_0]} \\
{\tilde{r}_0 = 0} & \quad & C_0 = Id \quad\quad\quad  C_{k+1} \in {[A_k]}{C_k}.
\end{eqnarray*}
Again, $[\tilde{r}]$ is evaluated by any method 0-2. Please
note that there is no need to inverse a matrix in the
doubleton representation when $[\tilde{r}]$ is evaluated
either by Method~0 or Method~2, so this approach is suitable
for the case where $A_k$ may be close to singular. In the
computer assisted proofs presented in this paper we use Method~0 to
represent $[\tilde{r}]$ because it is less computationally expensive
and, in our current setting, using the other methods have not improved the results.
This is puzzling, as it contradicts our experience with ODEs where Method~2 is preferable, and thus
it might be worthwhile to investigate this phenomena in some later study.

\subsection{Optimization exploiting the block structure of $D\Phi(x)$}
\label{subsec:optimization}

In our setting, $[A_k] = D\Phi(\bar{x}_k)$, where $\bar{x}_k$ is the $(p,n)$-f-set in the $k$-th step of integration.
As $[A_k]$ is $m \times m$ matrix, where $m$ is the size of $(p,n)$-representation ($m = p \cdot (n+2) + 1$), therefore, if
we decide to represent the error part $[\tilde{r}]$ in doubleton by interval box (Lohner Method 0), then
the matrix multiplications involving the matrix $[A_k]$ take the most of the execution time in one step of integration in the Lohner algorithm,
especially for large $n$ and/or $p$. From  Equations~(\ref{eq:phi-start}-\ref{eq:phi-end}),
we see that $D\Phi(x)$ has a nice block structure and contains a large number of zero entries, i.e. it is a sparse matrix.
This structure is well visible when we use the following index function $\mathbb{I}$:
\begin{eqnarray*}
\mathbb{I}(0,0) &=& 1, \\
\mathbb{I}(i,k) &=& 1 + (i-1) \cdot (n+1) + k, \quad 1 \le i \le p, 0 \le k \le n, \\
\mathbb{I}(i,n+1) &=& 1 + p \cdot (n+1) + i, \quad 1 \le i \le p.
\end{eqnarray*}
Under this index function, $[A_k]$ is of the form:
\begin{equation*}
[A_k] = \left(\begin{array}{cccc}
A_{11} & 0 & A_{13} & 0 \\
0 & Id & 0 & 0 \\
0 & 0 & 0 & 0
\end{array}\right),
\end{equation*}
where $A_{11}$ is $n+1 \times 1$ matrix (column vector), $A_{13}$ is $n+1 \times n+2$, and $Id$ is an identity matrix of size $(p-1) \cdot (n+1)$.
Therefore, we use this index function to define blocks for all matrices and vectors appearing in all methods discussed in Section~\ref{subsec:wrapping},
as, with such block representation of matrices and vectors, we can easily program the multiplication by a matrix or a vector
so that all the operations on any zero block are avoided. We will refer to this as the optimized algorithm.

If we have an arbitrary matrix $C$, then the cost of computing $[A_k] \cdot C$ by a standard algorithm for the matrix multiplication is of order $O(n^3 \cdot p^3)$
in both the scalar addition and multiplication operations (we remind, that $p$, $n$ are the parameters of $(p,n)$-representation).
In the case of the optimized algorithm, the block structure and sparseness of $[A_k]$ reduce the computational cost to
 $O(n^2 \cdot p^2)$ in scalar additions and $O(n^3)$ in scalar multiplications.

The computation times for the computer assisted proofs discussed in Section~\ref{sec:results} on
the $2.50$ GHz processor (see Section~\ref{sec:results} for a detailed specification) are presented in
Table~\ref{tab:times}. We see that the optimized algorithm is much faster than the direct multiplication, the
 speed up is evident especially for the larger $(p,n)$-representations.

\begin{table}[t]
\begin{center}
\caption{A comparison of the execution times for computer assisted proofs for standard and optimized matrix multiplication on 2.5GHz processor (no multi-threading).}
\label{tab:times}
\begin{tabular}{|l|c|c|}
\hline \textbf{Proof (p, n)} & \textbf{standard multiplication} & \textbf{optimized multiplication} \\
\hline Theorem~\ref{thm:mackey-glass-n6-stable} $(32, 4)$ & 45 seconds & 12 seconds \\
\hline Theorem~\ref{thm:mackey-glass-n8-stable} $(128, 4)$ & 133 minutes & 12 minutes \\
\hline
\end{tabular}
\end{center}
\end{table}

\section{Poincar\'e map for delay differential equations}
\label{sec:poinc-map}

\subsection{Definition of a Poincar\'e map}

We begin with the definition of the (transversal) section of the
semiflow $\varphi$ associated to \eqref{eq:dde}. First, we would like
to recall the  ODE setting where, for a flow
$\varphi : \R \times \R^m \to \R^m$, a (local) transversal section $S$
is usually defined as a (subset of) smooth manifold $\mathcal{M}$ of
codimension one satisfying the transversality condition:
\begin{equation}
\label{eq:ode-trans-cond}
\frac{d}{dt}\varphi(t, x_0)|_{t=0} = f(x_0) \not\in T_{x_0}\mathcal{M}, \quad x_0 \in S,
\end{equation}
where $T_x\mathcal{M}$ denotes the tangent bundle at $x$.
If $\mathcal{M}$ is a hyperplane
\begin{equation*}
\mathcal{M} = \left\{ x \in \R^m \ | \ s \cdot x = a \right\}
\end{equation*}
for some given normal vector $s \in R^m$, $\| s \| \neq 0$ and $a \in \R$,
then  condition~\eqref{eq:ode-trans-cond} becomes
\begin{equation*}
s \cdot f(x_0) \ne 0, \quad x_0 \in S.
\end{equation*}

We will use a similar approach in the context of the semiflow $\varphi$
associated to \eqref{eq:dde}. We will restrict ourselves to the linear
sections and we will use the fact that from Equation~\eqref{eq:dde-ivp} and
Lemma~\ref{lem:smoothing} it follows
\begin{equation*}
\frac{d}{dt}\varphi(t, x_0)|_{t=t_0} = \dot{x}_{t_0}
\end{equation*}
for any $t_0 \ge \tau$.
Moreover, $\dot{x}_{t_0}$ is of class $C^{n-1}$ wherever
$t_0 \ge n \cdot \tau$. This observation will be crucial for
the definition of a transversal section in the DDE context and,
later, for the rigorous computation of Poincar\'e maps.

\begin{definition}
\label{def:section}
Let $\varphi$ be the semiflow associated with the system \eqref{eq:dde}.
Let $n \in \mathbb{N}$, $s : C^n \to \R$ be a continuous  affine mapping,
i.e. $s(x)=l(x)+a$, where $l$ is a bounded linear functional and $a \in \mathbb{R}$.
We define \emph{a global \mbox{$C^n$-section}} as a hyperplane:
\begin{equation*}
\mathcal{S} = \{x \in C^n \ | \ s(x) = 0\}.
\end{equation*}
Any convex and bounded subset $S \subset \mathcal{S}$ is called
\emph{a local $C^n$-section} (or simply a section).

A section $S$ is said to be \emph{transversal} if
there exists a convex open set $W \subset C^n$ such that
\begin{equation*}
  W \cap \mathcal{S} = U, \qquad W=W_- \cup U \cup W_+,
\end{equation*}
where
\begin{equation*}
  W_- = \{x \in W \ | \  s(x)<0\}, \qquad W_+=\{x \in W \ | \ s(x)>0\}, \qquad \cl(U) = \cl(S),
\end{equation*}
satisfying the condition
\begin{equation}
\label{eq:sec-trans}
  l\left( \dot{x} \right) > 0, \quad \forall x \in W \cap C^{n+1}.
\end{equation}
We will refer to (\ref{eq:sec-trans}) as the \emph{transversality condition}.
\end{definition}

\begin{rem}
\label{rem:regularity}
Please note that the requirement $x \in W \cap C^{n+1}$ in \eqref{eq:sec-trans}
is essential to guarantee that $\dot{x}$ and thus $l(\dot{x})$ in \eqref{eq:sec-trans}
are well defined, as, for $x_0 \in C^{n+1}$ and $t > 0$, it might happen that $\varphi(t, x_0)$
is of class $C^k$, $k < n+1$ (the loss of regularity), but Lemma~\ref{lem:smoothing} states
that we only need to ,,long enough'' integrate the initial functions to get rid of this problem completely.
Those two phenomena are illustrated in the following example.

Let $x_0 \in W$ for $W$ as in
Definition~\ref{def:section}. In general, it happens
that $\varphi(\epsi, x_0) \not\in C^n$ for small $\epsi > 0$.
This may seem at first to contradict intuition from Lemma~\ref{lem:smoothing},
but in fact it is not. Consider the following r.h.s. of Equation~\eqref{eq:dde}:
\begin{equation*}
f(z_1, z_2) = z_1, \quad \forall z_1, z_2 \in \mathbb{R}
\end{equation*}
Let $x_0 \equiv 1$ be an initial function and let $x$ be a solution
to Equation~\eqref{eq:dde} with $x|_{[-1,0]} \equiv x_0$ and delay $\tau = 1$.
We see that $x(t)$ is $C^\infty$ on
$[-1, 0)$. However at $t=0$ we have
\begin{equation*}
0 = x'(0^-) \ne 1 = x'(0^+) = f\left(x(-1), x(0)\right),
\end{equation*}
so $x_t : [-1,0] \to \R$ is only of class $C^0$ for any
$t \in (0, 1)$. This is a very undesirable phenomena, but
 the solution will be smoothed
after a full delay, according to Lemma~\ref{lem:smoothing}.
As $x(t) = 1+t$ on $[0, 1]$, we have at $t = 1$:
\begin{equation*}
1 = x'(1^-) = x'(1^+) = f\left(x(0), x(1)\right).
\end{equation*}
One can show again that $x^{(2)}(1^-) \ne x^{(2)}(1^+)$, therefore $x$ is
of class $C^1$ on $(1, 2)$, and the smoothing of solutions goes on with the increasing $t$.

This shows for any $x_0 \in C^n$, if $\omega > n \cdot \tau$, then
we have ,,only'' $\varphi(\omega, x_0) \in C^n$ in a general case. On
the other hand, ,,long enough'' integration time $\omega$ can be used
to guarantee that every initial function $x \in C^0$ has a well defined
image in $C^n$ under mapping $\varphi(\omega, \cdot)$. This is essential in the
following construction of a Poincar\'e map for DDEs (Theorem~\ref{thm:tran-time}
and Definition~\ref{def:poincare-map}).
\end{rem}


\begin{theorem}
\label{thm:tran-time}
Assume $ n \in \mathbb{N}$, $V \subset C^0$. Let $S$ be a
local transversal $C^n$-section for \eqref{eq:dde}
and let $W$ be as in Definition~\ref{def:section}.
Let $\omega = (n+1) \cdot \tau$.

Assume that there exist $t_1,t_2 \in \R$, $\omega \le t_1 < t_2$, such that the following
conditions hold for all $x \in V$:
\begin{eqnarray}
\label{eq:darboux}
\varphi(t_1,x) \in W_- & \quad and \quad & \varphi(t_2,x) \in W_+ .
\end{eqnarray}

Then, for each $x_0 \in V$, there exists unique $t_{S}(x_0) \in (t_1, t_2)$
such that $\varphi\left(t_{S}(x_0), x_0\right) \in S$.
Also, $t_{S} : V \to [t_1, t_2]$ is continuous.
\end{theorem}
\textbf{Proof:}
Let $x_0 \in V$. By assumptions, $\varphi(t, x_0) = x_t \in W$
for $t \in [t_1, t_2]$ but also $x_t \in C^{n+1}$, by the assumption
on constants $\omega$, $t_1$, $t_2$ (by Lemma~\ref{lem:smoothing}).
So $s(\dot{x}_t)$ is well defined
and Condition~\eqref{eq:sec-trans} guarantees that
\begin{equation*}
\frac{d }{dt} s(\varphi(t,x)) = l\left(\frac{\partial}{\partial t}{\varphi}(t,x)\right)=l(\dot{x}_t) > 0, \qquad t \in [t_1,t_2].
\end{equation*}
Therefore, the function defined by $\hat{s}(t) := s(x_t)$ is continuous
and strictly increasing on $[t_1, t_2]$. Now, from
\eqref{eq:darboux}, it follows there exists unique
$t_0 \in (t_1, t_2)$ such that
$s\left(\varphi(t_0, x_0)\right) = c$. Together with
continuity of $\varphi$ (Lemma~\ref{lem:semiflow}), this
implies continuity of $t_{S}:V \to (t_1,t_2)$.
\qed

\begin{definition}
\label{def:poincare-map}
The same assumptions as in Theorem~\ref{thm:tran-time},
in particular assume $t_1 \ge \tau \cdot (n+1) = \omega$.
We define \emph{the transition map to the section $S$
after (at least) $\omega$} by
\begin{equation*}
P_{\ge \omega} : V \to S \subset C^n, \quad P_{\ge \omega}(x_0) := \varphi\left(t_{S}(x_0), x_0\right),
\end{equation*}
for $t_S$ defined as in Theorem~\ref{thm:tran-time}.
If  $V \subset S$, then the map $P_{\ge \omega}$
will be called \emph{the Poincar\'e return map on the section $S$
after $\omega$}.
\end{definition}

Finally, we state the last and the most important theorem that
will allow us to apply  topological  fixed point theorems to $P_{\ge \omega}$.

\begin{theorem}
\label{thm:poincare-comp}
 Consider Poincar\'e map (after $\omega$) $P_{\ge\omega}: S \supset V \to S$ for some section $S$ under
the same assumptions as in Theorem~\ref{thm:tran-time}, especially assume $(n+1) \cdot \tau \le \omega \le t_S(V) \in [t_1, t_2]$.

Assume additionally that $\varphi([\omega,t_2],V)$ is bounded in $C^n$.

Then the map $P_{\ge\omega}$ is
continuous and compact in $C^n$, i.e. if $K \subset V$ is bounded,
then $\cl (P_{\ge\omega}(K))$ is compact in $C^n$.
\end{theorem}
\textbf{Proof:}

By Theorem~\ref{thm:tran-time}, $P_{\ge\omega}$ is well defined
for any $x_0 \in V$ since $t_1 \ge \omega \ge (n+1) \cdot \tau$ and
$P_{\ge\omega}(x_0) \in C^{n+1}$.

The continuity follows immediately from the continuity of
$t_{S}$ (Theorem~\ref{thm:tran-time}) and $\varphi$ (Lemma~\ref{lem:semiflow}).

Let $D = P_{\ge\omega}(V)$. From our assumptions, it follows that $D$ is bounded in $C^n$.
A known consequence of the Arzela-Ascoli Theorem is that, if
$D \subset C^{n}$ is closed and bounded, $x \in D$, $x^{(n+1)}$ exists, and there is $M$ such that
$\sup_{t} \left|x^{(n+1)}(t)\right| \leq M$ for all $x \in D$, then $D$ is compact (in $C^{n}$-norm).
Therefore, to finish the proof, it is enough to
 show that there is a uniform bound on $P_{\ge\omega}(x)^{(n+1)}$. For this, it is sufficient to have
 a uniform bound on $(\varphi(t,x))^{(n+1)}$ for $t \in [\omega,\sup_{x \in V} t_S]$. The existence of this bound
 follows from boundedness of derivatives up to order $n$ and formula \eqref{eq:F[k]=x[k]}.

\qed

The restriction on the transition time may seem a bit
unnatural since each solution becomes
 $C^\infty$ eventually,
as discussed in Remark~\ref{rem:regularity}.
In fact, it should be possible to work directly with the
solutions on the $C^n$ solution manifold $M^n$ (i.e. $M^n \subset C^n$ and $\varphi(t, M^n) \subset M^n$ for all $t \ge 0$).
When we restrict the flow to the solutions manifold $M^n$, then we do not need to demand
that the transition time to the section is bigger than $\omega = (n+1) \cdot \tau$.
Instead, to obtain the compactness, we need to shift the set forward only by one full delay.
Therefore we obtain the following theoretical result:
\begin{theorem}
\label{thm:poincare-comp-Mn}
 Consider Poincar\'e map (after $\tau$)
$P_{\ge\tau}: S \supset V \to S$ for some section $S$,
where $V \subset M^n$. Let $t_1$ and $t_2$ be like
in Theorem~\ref{thm:tran-time}.

Assume that $\varphi([0,t_2],V)$ is bounded in $C^n$.

Then the map $P_{\ge\tau}: V \to S \cap M^n$ is
continuous and compact in $C^n$, i.e. if
$K \subset V$ is bounded, then
$\cl (P_{\ge\omega}(K))$ is compact in $C^n$.
\end{theorem}

At the present stage of the development of our algorithm, we do  not have  the
constructive parametrisation of the manifold $M^n$, therefore we need to use the ,,long enough''
integration time $\omega$ in the rigorous numerical computations.

\subsection{Rigorous computation of Poincar\'e maps}
\label{sec:poincare-algorithm}

The restriction of the integration procedure $I_h$
(Section~\ref{sec:integrator}) to the fixed-size step $h =
\frac{1}{p}$ is a serious obstacle when we consider computation of
Poincar\'e map $P_{\ge\omega} : V \to S$. Obviously, if we assume for
simplicity that $\omega = q \cdot {\tau \over p}$, $q \in \mathbb{N}$
and $t_S(V) \subset \omega + [\epsi_1, \epsi_2]$ with
$0 < \epsi_1 < \epsi_2< {\tau \over p}$, then we have to find a method
to compute image of the set after small time $\epsi \in [\epsi_1, \epsi_2]$. The
definition of the \mbox{$(p,n)$-representation} together with
Equation~\eqref{eq:representation-taylor-interpretation} give a hint
how to compute the value of the function (and the derivatives up to the order $n$) for
some intermediate time $0 < \epsi < h$. But again, we face yet another
obstacle, as computing the $(p,n)$-f-set representing $\varphi(\epsi, x_0)$
for all initial functions $x_0$ in some given $(p,n)$-f-set turns out to
be impossible. It can be seen from the very same example as in
Remark~\ref{rem:regularity}. In the example, $x_\epsi$ would be only
$C^0$ at $t = -\epsi$. So if $\epsi$ is not a multiple of $h$, then, for
any $n > 0$, there is no $(p,n)$-representation of $x_\epsi$, unless we
restrict the computations to the set $C^n_p \cap C^{n+1}$
(or to the solutions manifold $M^n$). This is again a reason for
an appearance of the ,,long enough'' integration time in
Definition~\ref{def:poincare-map}.

This discussion motivates the following definition and lemma.
\begin{definition}
 Let $\bar{x}_0$ be a $(p,n)$-representation, and let
$\bar{x}_h = I_h(\bar{x}_0)$ and $c_{\bar{x}_0}$ be as in Definition~\ref{def:cksymb}.
For $\epsi \in [0, h]$ we define $(p,n)$-f-set $\bar{x}_\epsi$ by the following $(p,n)$-representation:
\begin{eqnarray*}
\repr{i}{k}{\epsi} & := & c^{i,[k]}_{\bar{x}_0}(\epsi), \quad i \in \{1,..,p\}, k \in \{0,..,n+1\}, 	       \\
\repr{i}{n+1}{\epsi} & := & hull\left( \repr{i}{n+1}{0}, \repr{i}{n+1}{h} \right), \quad i \in \{1,..,p\},  \\
\repr{0}{0}{\epsi} & := & \sum_{k=0}^{n+1} \repr{1}{k}{h} \cdot \varepsilon^k.
\end{eqnarray*}
For a given $\bar{x}_0$ we denote:
\begin{equation*}
I_\epsi(\bar{x}_0) = \bar{x}_\epsi
\end{equation*}
Function $I_\epsi$ will be called \emph{the shift by $\epsi$} or \emph{the \epsi-step integrator}.
\end{definition}

\begin{rem}
$I_\epsi(\bar{x}_0)$ is constructed in such a way that it contains
all solutions to \eqref{eq:dde-ivp} for initial functions
$x_0 \in Supp(\bar{x}_0) \cap M^{n+1}$ after time $t = \epsi$.
\end{rem}

\begin{theorem}
\label{thm:epsi-step}
Assume that $\epsi \in [0, h]$, $\bar{x}_0$ is a $(p,n)$-f-set. Let define
\begin{equation*}
\bar{x}_j = I_h(\bar{x}_{j-1}), \quad j \in \mathbb{N}.
\end{equation*}

If $n \cdot p \le q \in \mathbb{N}$, then
\begin{equation*}
\varphi(q \cdot h + \epsi, x) \in I_\epsi( \bar{x}_{q} )
\end{equation*}
for all $x \in \bar{x}_0$.
\end{theorem}
\textbf{Proof:} Since $q \ge n \cdot p$ then $q \cdot h \ge n \cdot \tau$
and the proof follows from Lemma~\ref{lem:smoothing}, Lemma~\ref{lem:c-k}
and Definition~\ref{def:representation}.
\qed

Now, the application of $I_h$ and $I_\epsi$ to compute $P_{\ge\omega}$ is straightforward.

\vskip 1em \noindent
{\large\textbf{Program $P_{\ge\omega}$} }
\vskip 0.5em \noindent

\noindent \textbf{Input: }
\begin{enumerate}
\item
a section $S$;

\item
a $(p,n)$-f-set $\bar{x}_0 \subset S$;

\item
$\omega = (n+1) \cdot \tau$;
\end{enumerate}

\noindent \textbf{Output: }
\begin{enumerate}
\item
$q \in \mathbb{N}$, $0 < \epsi_1 < \epsi_2 < {\tau \over p}$
such that $t_S(\bar{x}_0) \subset q \cdot {\tau \over p} + [\epsi_1, \epsi_2]$
for $\omega \le q \cdot {\tau \over p}$;

\item
$(p,n)$-f-set $\bar{x}_{\ge\omega}$ such that $P_{\ge\omega}(x_0) \in \bar{x}_{\ge\omega}$ for all $x_0 \in \bar{x}_0$;
\end{enumerate}

\noindent \textbf{Algorithm: }
\begin{enumerate}
\item do at least $(n+1) \cdot p$ iterations of $I_h$
to guarantee the $C^{n+1}$  regularity  of solutions for all initial functions
(so the map $P_{\ge (n+1) \cdot \tau}$ is well
defined and compact).

\item find $q \ge (n+1) \cdot p$ and $\epsi_1$, $\epsi_2 < h$
(for example by the binary search algorithm) such that for
the assumptions of Theorem~\ref{thm:tran-time} are guaranteed for
section $S$, $t_1 = q \cdot h + \epsi_1$, $t_2 = q \cdot h + \epsi_2$ and set $W$ defined by
\begin{equation*}
W := C^n \cap \left( I_{[\epsi_1, \epsi_2]} \circ I_h^q (\bar{x}_0) \right).
\end{equation*}

\item  By assumptions and by Theorems~\ref{thm:tran-time}~and~\ref{thm:epsi-step},
we know that we have $P_{t_1}(x_0) \in W \cap S$ for each $x_0 \in V$.
Moreover, by Theorem~\ref{thm:poincare-comp}, the map $P_{t_1}$ is compact (in $C^n$).

\end{enumerate}

\vspace{0.5cm}

Please note, that the operator $I_{[\epsi_1, \epsi_2]}$ should be interpreted as
computation of the sum $\bigcup_{\epsi \in [\epsi_1, \epsi_2]} I_\epsi(\cdot)$
or as any reasonable bound on this sum. In our program we just
evaluate $I_\epsi$ with $\epsi=[\epsi_1, \epsi_2]$.

\begin{rem}[Controlling the wrapping effect for $I_\epsi$]
We can use the decomposition of $I_\epsi$ into $\Phi_\epsi$ and
$R_\epsi$ such that the Lohner algorithm can again be used
in the last step of the integration as described in Section~\ref{subsec:wrapping}.
We skip the details and refer to the source code documentation of
the library available at \cite{szczelina-phd}.
\end{rem}

 Now, the question arises: how to represent the section $S$
in a manner suitable for computation of the program $P_\ge \omega$?

\subsection{$(p,n)$-sections}
Since we are using the $(p,n)$-representations to describe functions in $C^n$,
it is advisable to define sections in such a way that it would be easy
to rigorously check whether $x \in S$ for all functions represented by
a given $(p,n)$-f-set. The  straightforward way is to require
$l$ in the definition of $S$ to depend only on representation coefficients
$\repr{i}{k}{}$.

\begin{definition}
\label{def:pn-sec}
Let $l_{i,k}\in \R$ for $(i,k) \in \mathcal{C}=\{1,\dots,p\} \times \{0,\dots,n\} \cup \{(0,0)\}$. We assume that at least one $l_{i,k}$ is not equal to zero.
Let $a \in \R$ be given.
 For $x \in C^n$ we define a linear continuous map $l : C^n \to \R$ by
\begin{equation}
\label{eq:pnsection}
l(x) = \sum_{(i,k) \in \mathcal{C}} l_{i,k} x^{[k]}(-i \cdot h).
\end{equation}
The section $S=\{x \in C^n \ | \ l(x) - a = 0\}$  is called \emph{a $(p,n)$-section}.
\end{definition}

\subsection{Choosing an optimal section}
\label{sec:good-section}

\begin{figure}
    \centering{
        \includegraphics[width=50mm]{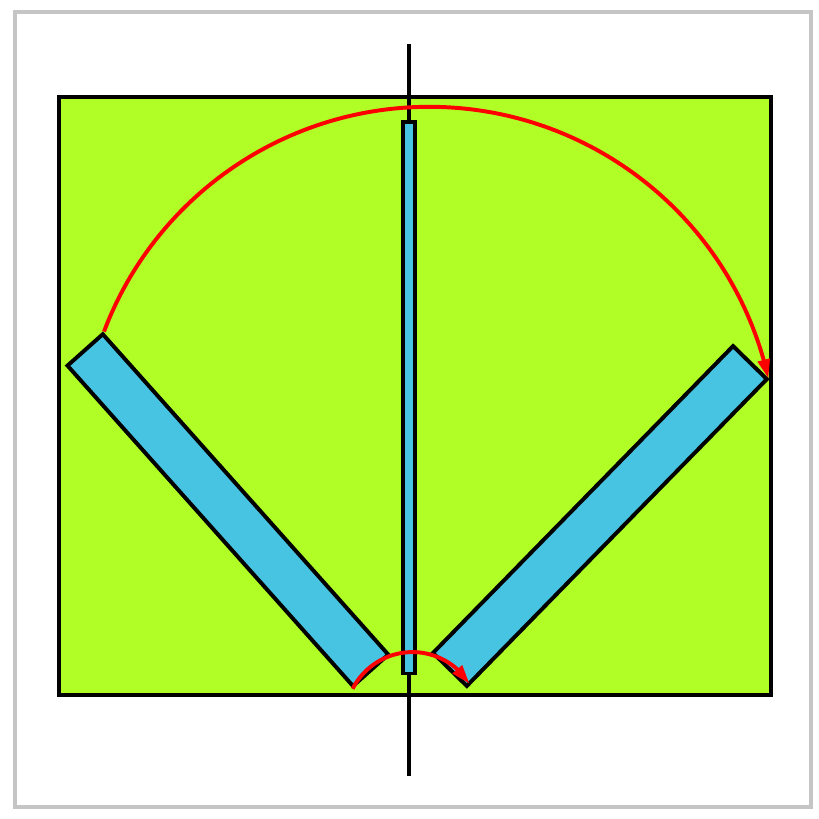}
        \includegraphics[width=50mm]{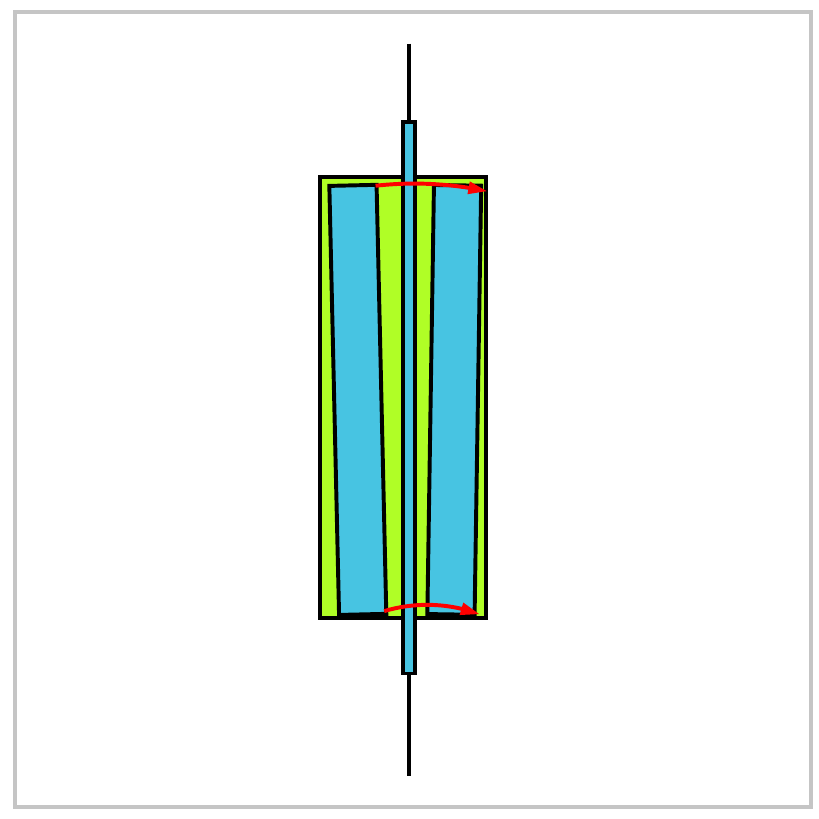}
    }

    \caption{\label{fig:section-choose}
	\textbf{Left: }an illustration of the problem with big difference in transition
     time $t_S$ for poorly chosen section $S$. If $t_S \in q \cdot {\tau \over p} + [\epsi_1, \epsi_2]$
     and $|\epsi_1 - \epsi_2|$ is large, then $I_{[\epsi_1, \epsi_2]}$ produces estimates on
     solutions distant from the section (blue rectangles) so the interval enclosure $W$ of all
     solutions tend to be very large (green rectangle). \textbf{Right: }if the section is chosen
     carefully, then all the
     solution obtained from $I_{[\epsi_1, \epsi_2]}$ are close to the section, so the set $W$ is
     small.
    }
\end{figure}

Numerical experiments with the rigorous integrator have
shown that the choice of a good section is a key factor
to obtain sufficiently good bounds on images of the Poincar\'e map
to be used in computer assisted proofs reported in Section~\ref{sec:results}.
The section has to be chosen so that  the diameter of the bounds on transition
time $t_S$ should be as small as possible, see Figure~\ref{fig:section-choose}.

We will discuss the problem of choosing optimal section in the ODEs case.
Later, in Section~\ref{sec:results}, we will apply a heuristic procedure based on this
discussion to obtain a good candidate for an optimal section in the DDEs setting.

Let us consider an ODE of the form:
\begin{equation}
 \label{eq:ode}
  x'=f(x), \quad f\in C^1, \quad x\in \R^n.
\end{equation}
Let $x_0$ be a periodic orbit of period $T$ of the flow $\varphi$
induced by \eqref{eq:ode}. Then, $f(x_0)$ is a right eigenvector of the
matrix $M = \frac{\partial \varphi}{\partial x}(T,x_0)$ with eigenvalue
$\lambda=1$. Let $l$ be a row vector which is a left eigenvector of
$M$ corresponding to $\lambda=1$.  Let us assume, that the periodic orbit passing
through $x_0$ is hyperbolic. In such a case, the left and right eigenvectors
corresponding to the eigenvalue $\lambda=1$ are uniquely defined up to
a multiplier and we have
\begin{equation*}
  l \cdot f(x_0) \neq 0.
\end{equation*}
We normalize $l$ so that
\begin{equation*}
  l \cdot f(x_0)=1.
\end{equation*}

For any given row vector $v \in \mathbb{R}^n$ let us consider a
section $S_v=\{ x \ | \ v \cdot x - v \cdot x_0=0 \}$.
We define
\begin{equation*}
  v^{\bot}=\{x \in \R^n \ | \  v \cdot x =0\},
\end{equation*}
hence $v^\bot$ is the tangent space to the section $S_v$.

Under the above assumptions, we have the following lemma.
\begin{lemma}
\label{lem:good-section}
If $v \cdot f(x_0) \neq 0$, then $S_v$ is locally transversal and
\begin{equation}
  \label{eq:d-tran-time}
  \frac{\partial t_{S_v}}{\partial x}(x_0)=-\frac{v \cdot \frac{\partial \varphi}{\partial x}(T,x_0)}{v \cdot f(x_0)},
\end{equation}
where $t_{S_v}$ is the transition time to the section $S_v$,
defined in some neighborhood of $x_0$.

Moreover,
\begin{equation}
  \label{eq:dtran=0}
  \frac{\partial t_{S_v}}{\partial x}(x_0) \cdot b=0, \quad \mbox{for $b \in v^{\bot}$},
\end{equation}
iff $v = \alpha l$  for some $\alpha \neq 0$
\end{lemma}
\textbf{Proof:}
The transition time to section $S_v$ is defined by the following implicit equation
\begin{equation*}
  v \cdot \varphi(t_{S_v}(x),x) - v \cdot x_0=0.
\end{equation*}
From this, we immediately obtain \eqref{eq:d-tran-time}.

The second assertion is obtained as follows. At first, assume \eqref{eq:dtran=0}. We have
\begin{eqnarray*}
   \frac{\partial t_{S_v}}{\partial x}(x_0) \cdot f(x_0) = -\frac{v \cdot \frac{\partial \varphi}{\partial x}(T,x_0)}{v \cdot f(x_0)} \cdot f(x_0) =  \\
   - \frac{1}{v \cdot f(x_0)} v \cdot \left( \frac{\partial \varphi}{\partial x}(T,x_0) \cdot f(x_0) \right) = \ - \frac{1}{v \cdot f(x_0)} \cdot (v \cdot f(x_0)) = -1.
\end{eqnarray*}
Therefore $v$ is proportional to $l$.

The other direction of the second assertion is obvious.
\qed

In simple words, Lemma~\ref{lem:good-section} states that choosing the
left eigenvector of the monodromy matrix $\frac{\partial \varphi}{\partial x}(T,x_0)$
gives a section such that the return
time to this section is constant in the first order
approximation.

\section{The existence of periodic orbits in Mackey-Glass equation}
\label{sec:results}

The Mackey-Glass system \eqref{eq:mackey-glass} is one of the best known
delay differential equations. The original work of Mackey and Glass \cite{mackey-glass} spawned
wide attention, being cited by many papers with a broad spectrum of topics:
from theoretical mathematical works to neural networks and electrical engineering.
 Numerical experiments show that, as either
parameter $\tau$ \cite{mackey-glass} or $n$ \cite{mackey-glass-scholarpedia}
is increased, the system undergoes a series of period doubling bifurcations
and they lead to the creation of an apparent chaotic attractor.

In this section, we present computer assisted proofs of the existence of attracting periodic
orbits in Mackey-Glass system \eqref{eq:mackey-glass}. We use the classical values of parameters:
$\tau = 2$, $\beta = 2$ and $\gamma = 1$
and we investigate the existence of periodic orbits with  $n = 6$ (before the first period doubling) and $n = 8$ (after
the first period doubling) \cite{mackey-glass-scholarpedia}. We would like to stress, that we are not proving
that these orbits are attracting. This would require some $C^1$-estimates for the Poincar\'e map
defined by \eqref{eq:mackey-glass}.

\subsection{Outline of the method for proving periodic orbits}

The scheme of a computer assisted proof of a periodic orbit consists of several steps:
\begin{enumerate}
\item \label{idea:finite} find a good, finite representation of
bounded sets in the phase space $C^k$ (or in other suitable function space),

\item \label{idea:initial} choose suitable  section $S$ and some
\emph{a priori} initial set $V$ on the section,

\item \label{idea:poincare} compute image of $V$ by Poincar\'e  map $P_{\ge \omega}$ on section $S$,

\item \label{idea:fixed-point} prove that the map $P_{\ge\omega}$, the set $V$, and
the set $W := P_{\ge\omega}(V)$ all satisfy assumptions of some fixed point theorem
so that it implies the existence of a fixed point for $P_{\ge\omega}$ in
$V$. This gives rise to the periodic orbit in Equation~\eqref{eq:dde}.
\end{enumerate}
To this point, we have presented ingredients needed in
steps~\ref{idea:finite}~and~\ref{idea:poincare}.
In Step~\ref{idea:fixed-point}, we will
use the Schauder Fixed Point Theorem \cite{Schauder,fixed-points-book}:
\begin{theorem}[Schauder Fixed Point Theorem]
\label{thm:schauder}
Let $X$ be a Banach space, let $V \subset X$ be non-empty, convex, bounded set
and let $P : V \to X$ be continuous mapping such that $P(V) \subset K \subset V$
and $K$ is compact. Then the map $P$ has a fixed point in $V$.
\end{theorem}
Theorem~\ref{thm:schauder} is suitable for proving the existence of periodic orbits
for which there is a numerical evidence  that they
are attracting.  The unstable periodic orbits can be treated
by adopting the covering relations approach from \cite{cov-rel}, which may
be applied in the context of infinite dimensional phase-space
(for such an adaptation in the context of dissipative PDEs see
\cite{zgliczynski-pde-3}).

In Section \ref{sec:good-section} we have presented some theoretical background on the
selection of a suitable section that is the foundation of Step~\ref{idea:initial}.
Now, we would like to put more emphasis on technical details,
as the procedure in Step~\ref{idea:initial} introduces some difficulties due to
large size of the data defining  $(p,n)$-representations. In the proofs we use $(32,4)$- and
$(128,4)$-f-sets with representation sizes $m = 193$ and $m = 769$, respectively.
Thus we are not able to simply ,,guess'' good coordinates or refine them ,,by hand'' -
we need an \emph{automated} way to do that.

The following discussion is a bit technical and involves some heuristics, thus it is probably 
relevant only for people interested in implementing their own version of the software.
Those interested only in the actual proofs of the existence of periodic orbits should move
to Section~\ref{sec:proofs}.

\subsection{Finding suitable section and good initial set for a computer-assisted proof}
\label{sec:selection-initial}
Here we give an outline for the selection of a good initial data for the proof of the existence of 
an apparently attracting orbit. It consists of the following steps:
\begin{enumerate}
\item\label{find-1} find a good numerical approximation
$x_0$ of a periodic solution to equation~\eqref{eq:dde},

\item\label{find-2} find a good section $S$ - by this we mean the difference
between transition times $t_1$ and $t_2$
(as defined in Theorem~\ref{thm:tran-time})
is as small as possible (in the vicinity of $x_0$),

\item\label{find-3} choose a good coordinate frame in $S$ for
the initial $(p,n)$-representation, then choose the $(p,n)$-f-set
$V \subset S$, such that  $x_0 \in V$ and  $P_{\ge\omega}(V) \subset V$.
\end{enumerate}
We will now describe shortly how each of the above steps was implemented.
In this description, we refer to non-rigorous computations, that is algorithms:
$\hat{I}_h$ defined as in Section~\ref{sec:integrator},
and $\hat{I}_\epsi$ defined as in Section~\ref{sec:poinc-map}, but
with the remainder terms depending on the rough-enclosure ignored
and explicitly set to $0$. Using non-rigorous integrators $\hat{I}_h$ and $\hat{I}_\epsi$, we
construct a finite-dimensional semiflow $\hat{\varphi}$ that approximates
$\varphi$ by
\begin{equation*}
\hat{\varphi}(t, x) := I_\epsi \circ I_h^q(x) \qquad t = q \cdot h + \epsi, \ \ \epsi < h.
\end{equation*}
Now, the procedure for finding good initial conditions can be described as follows:

\noindent\textbf{Step~\ref{find-1}.}
Since we are looking for an attracting orbit we start by non-rigorously
integrating forward in time an initial
function $\hat{x}_0 \equiv 1.1$ for some arbitrary, long time $T_{iter}$, 
until we see that $\hat{x}_{T_{iter}}$ approach the apparently stable periodic orbit. 
Then, we refine $\hat{x}_{T_{iter}}$ by the Newton algorithm
applied to $x \mapsto \hat{P}_{\ge\omega}(x) - x$, where the map $\hat{P}_{\ge\omega}$
is a non-rigorous version of $P_{\ge\omega}$ defined as a first return map for semiflow $\hat{\varphi}$
to a simple section
$S = \{x \ | \ x(0) = \hat{x}_{T_{iter}}(0)\}$. The output of this step is
a numerical candidate for the periodic solution $x_0$, given by its $(p,n)$-representation $\bar{x}_0$
such that $\bar{x}_0$ and $P_{\ge\omega}(\bar{x}_0)$ are close.

\noindent\textbf{Step~\ref{find-2}.}
This is an essential step, as numerical experiments with the
rigorous integrator have shown that the choice of a good
section is the key factor to obtain tight bounds
on the image of the Poincar\'e map. We use an observation from Section~\ref{sec:good-section}
and we find the left eigenvector $\hat{l}$ of the matrix
$\frac{\partial \hat{\varphi}}{\partial x}(T,x_0)$ corresponding to eigenvalue $1$,
where $T$ is an apparent period of the approximate periodic orbit for the non-rigorous semiflow $\hat{\varphi}$.

 Please note that $\hat{l}$ might be considered a $(p,n)$-representation with
remainder part set to $0$, therefore we can define a $(p,n)$-section by
\begin{equation}
  \hat{l}\cdot x - \hat{l}\cdot \bar{x}_0=0, \label{eq:sec-l-vect}
\end{equation}
where the dot product is computed using the coordinates of $(p,n)$-representation, i.e.  in the vector space $\R^m$, where $m$ is
the size of a $(p,n)$-representation, $m = (n+2) \cdot p + 1$.

\noindent\textbf{Step~\ref{find-3}.}
Having a good candidate for the section $S$ (defined by (\ref{eq:sec-l-vect})), we need to introduce the coordinates on it.
For this, we create the following matrix:
\begin{equation*}
  A := \left( \ \hat{l}^T \ \ \vrule \ \begin{array}{c} 0 \quad ... \quad 0 \\ Id_{m-1 \times m-1} \end{array} \right),
\end{equation*}
Now, let $C$ denote the matrix obtained after orthonormalization of columns  of $A$. Please note,
that matrix $C$ acts on the variables corresponding to the remainder terms as an identity.
This follows from the fact that $\hat{l}_{i,[n+1]}\equiv 0$.
It is easy to see that all $(p,n)$-representations that lie on the section $S$ are given by:
\begin{equation*}
  y \mapsto \bar{x}_0 + C y,
\end{equation*}
for all $y$ such that $\pi_1 y = 0$.

Now, on section $S$, using the coordinates defined by the matrix $C$, we define a candidate set $[V]$, in a form of $(p,n)$-f-set in a following manner.
Let $[r] \subset \R^{m-p}$ (these correspond to variables  $x^{i,[k]}$ for $k \leq n$) and $[B] \subset \R^{p}$ (these are bounds for $x^{i,[n+1]}$ - the remainders)
be two interval boxes centered at $0$ such that  $diam(\pi_1 [r]) = 0$. We put $[\bar{r}_0] := [r] \times [B]$
and we define $(p,n)$-f-set $[V]$ by:
\begin{equation*}
  [V] := \bar{x}_0 + C \cdot [\bar{r}_0] .
\end{equation*}

Diameters of $\pi_i [\bar{r}_0]$ for $i \geq 2$ are selected
experimentally to follow some exponential law in parameter $k$
(i.e. $diam(\pi_{\mathbb{I}(i, k)}[\bar{r}_0]) \approx a^k$ for $1 \le i \le n$), as
the periodic solutions to Equation~\eqref{eq:mackey-glass} are at least of class $C^\infty$
and, if $x(t) > -1$ for all $t$, then they should be analytic \cite{mallet-parret-analytic, nussbaum-analytic}.
The remainder $[B]$ is chosen initially such that $diam\left(B\right) \gg diam\left([r_0]\right)$.
Therefore the initial selection of $[\bar{r}_0]$ may not be good enough to
satisfy assumptions of Theorem~\ref{thm:schauder} right away. As the dynamics of the system is strongly contracting, we hope
to obtain a good initial condition by the following iteration. We start with $[V]_0 = [V]$ and we compute
$[V]_{i+1} = P_{\ge\omega}([V]_i) \cap [V]_i$, until the condition $P_{\ge\omega}([V]_i) \subset [V]_i$
is eventually meet at some $i_{stop}$. Then the initial set for the computer assisted proof is $[\tilde{V}] = V_{i_{stop}}$.
Both initial sets that are used in computer assisted proofs in this paper were generated with such a procedure (see source codes).

Observe that we are not very careful in the choice of \emph{coordinates \underline{on} the section} - we simply
choose some basis orthonormal to the normal vector $\hat{l}$ of the section hyperplane.
Definitely better choice would be to use approximate eigenvectors of the Poincar\'e map,
but in the case of strongly attracting periodic orbits it is enough to choose a good section.
Observe also, that the orthonormal matrix is easy to invert rigorously, which is an important step
in comparison of the initial set and its image by the Poincar\'e map.

\subsection{Attracting periodic orbits in Mackey-Glass equation for $n = 6$ and $n = 8$}
\label{sec:proofs}

In this section we present two theorems  about the existence of periodic orbits in
Mackey-Glass equation. As they depend heavily on the estimates obtained from the rigorous
numerical computations, we would like to discuss first the textual presentation of numbers
used in this section and how they are related to the input / output
values used in rigorous computations.

In the rigorous numerics we use intervals with ends being representable computer numbers.
The representable numbers are implemented as \verb|binary32| or \verb|binary64| data types defined in
IEEE Standard for Floating-Point Arithmetic (IEEE 754) \cite{ieee-754}, so that they are
stored (roughly speaking) as $s \cdot m \cdot 2^e$, where $s$ is the sign bit, $m$ is the mantissa
and $e$ the exponent. Such a representation means that most numbers with a finite representation in the decimal base
are not representable (e.g. number $0.3333$). In this paper, for better readability,
we are going to use the decimal representation of numbers with the fixed precision
(usually $4$ decimal places), so we have rewritten computer programs to handle those values rigorously.
For example, if we write in the text that $a = 0.3333$ then we put the following rigorous operation
in the code:
\begin{equation*}
a = [3333,3333] \div [10000, 10000].
\end{equation*}
That is, all numbers presented here in theorems and/or proofs should be
regarded by the reader as the real, rigorous values, even if they are not representable in the sense of  IEEE 754 standard.

In the proofs we refer to computer programs \verb|mg_stable_n6| and \verb|mg_stable_n8|.
Their source codes, together with instructions on the compilation process,
can be downloaded from \cite{szczelina-www-mackeyglass}.
The codes were tested on a laptop with
Intel\textsuperscript{\textregistered}
Core\textsuperscript{\texttrademark}
I7-2860QM CPU (2.50 GHz),
16 GB RAM under 64-bit Linux operating system
(Ubuntu 12.04 LTS) and
C/C++ compiler \verb|gcc| version 4.6.3.

\subsubsection{Case $n=6$}
Our first result is for the periodic orbit for the parameter value before the first period doubling bifurcation.

With $n=6$, numerical experiments clearly show that the minimal period of the periodic orbit is around $5.58$.
In our proof however, due to the problem with the loss of  the regularity at the grid points, thus the need to use
the ,,long enough'' transition time, we consider the second return to the section.

Numerical experiments indicate that the orbit is attracting with the most significant
eigenvalues of the map $P_{\ge \omega}$ (again, this is the second return to the
Poincar\'e section) estimated to be:
{\begin{center}
\begin{tabular}{|c|c|c|c|c|c|c|c|c|c|}
\hline
$Re \lambda$ & -0.0437 & -0.0437 & 0.0030 &  0.0030 & -0.0028 & 0.0019 & -0.0003 & -0.0003 & 0.0005 \\
\hline
$Im \lambda$ &  0.0793 & -0.0793 & 0.0097 & -0.0097 &         &        &  0.0018 & -0.0018 &        \\
\hline
$ |\lambda$| &  0.0905 &  0.0905 & 0.0102 &  0.0102 &  0.0028 & 0.0019 &  0.0019 &  0.0019 & 0.0005 \\
\hline
\end{tabular}
\end{center}
}
Therefore, the contraction appears to be quite strong, so the
choice of good coordinates on the section appears to be not important.

We obtained the following theorem.
\begin{theorem}
\label{thm:mackey-glass-n6-stable}
There exists a $T$-periodic solution $x$ with period $T \in \left[10.9671, 10.9673\right]$
to Equation~\eqref{eq:mackey-glass} for parameters $\gamma = 2$, $\alpha = 1$, $\tau = 2$ and $n = 6$.
Moreover
\begin{eqnarray*}
    \left\| \hat{x} - x\right\|_{C^0}  &\le& 0.02 \\
    \left\| \hat{x} - x\right\|_{C^1}  &\le& 0.05 \\
    \left\| \hat{x} - x\right\|_{C^2}  &\le& 0.08 \\
    \left\| \hat{x} - x\right\|_{C^3}  &\le& 0.13 \\
    \left\| \hat{x} - x\right\|_{C^4}  &\le& 0.18
\end{eqnarray*}
for $\hat{x}$ defined by
\begin{equation*}
\renewcommand*{\arraystretch}{1.5}    \begin{array}{rclclc}
\hat{x}(t) & = & 0.9773 & - & &  \\
           & - & 0.0031 \cdot \cos\left( {2\pi \over T} \cdot 2 \cdot t \right) & + & 0.2398 \cdot \sin\left( {2\pi \over T} \cdot 2 \cdot t \right) & + \\
           & + & 0.0165 \cdot \cos\left( {2\pi \over T} \cdot 4 \cdot t \right) & - & 0.0043 \cdot \sin\left( {2\pi \over T} \cdot 4 \cdot t \right) & + \\
           & + & 0.0102 \cdot \cos\left( {2\pi \over T} \cdot 6 \cdot t \right) & - & 0.0011 \cdot \sin\left( {2\pi \over T} \cdot 6 \cdot t \right) & - \\
           & - & 0.0007 \cdot \cos\left( {2\pi \over T} \cdot 8 \cdot t \right) & + & 0.0014 \cdot \sin\left( {2\pi \over T} \cdot 8 \cdot t \right) &
\end{array}
\end{equation*}
\end{theorem}
\textbf{Proof:} Verification of assumptions of the Schauder theorem is done
with the computer assistance in the program \verb|mg_stable_n6|.
It uses the $(p,n)$-representation of the phase-space with $p = 32$ and $n = 4$.
Initial $(p,n)$-f-set $[\bar{x}_0]$ is provided directly in the source code and
it was selected with procedure described in Section~\ref{sec:selection-initial}.
For the  map $P_{\ge\omega}$, $\omega = (n+1) \cdot \tau = 10$, which represents the second return
to the section $S$, we obtained:
\begin{eqnarray*}
    t_S                         & \in & \omega + \left[\varepsilon\right] \ \subset \ \left[10.9671,10.9673\right],    \\
    \omega                      &  =  & q \cdot \frac{\tau}{p}         \     =    \ 175 \cdot \frac{\tau}{p}, \\
    \left[\varepsilon\right]    &  =  & \left[0.02960307544, 0.02971816895\right],
\end{eqnarray*}
which guarantees the $C^{n+1}$-regularity of the solutions
and the compactness of the map $P_{\ge\omega}$ (in $C^k$ norm for $k \le n$).
The inclusion condition $P_{\ge \omega}\left([\bar{x}_0]\right) \subset [\bar{x}_0]$ of the Schauder Fixed Point Theorem
is checked rigorously, see output of the program for details.  Together, these two
facts guarantee the non-emptiness of $Supp^{(n+1)}([\bar{x}_0])$.
The transversality is guaranteed with $l(\dot{x}) \ge 0.2828$ for all $x \in C^{n+1} \cap P(\bar{x}_0)$.
The distance in $C^{0}$ norm is rigorously estimated to
\begin{equation*}
\|\hat{x}-x\|_{C^0} \leq 0.01902867681.
\end{equation*}
Similarly, we have verified the other norms, see output of the program.
\qed

The execution of the program realizing this proof took around 12 seconds on $2.50$ GHz machine.

The diameter of the estimation for period $T$ (also for the last step $[\epsi_1, \epsi_2]$)
obtained from the computer-assisted proof is close to $1.15 \cdot 10^{-4}$.

\begin{figure}
\centering{
\includegraphics[width=125mm]{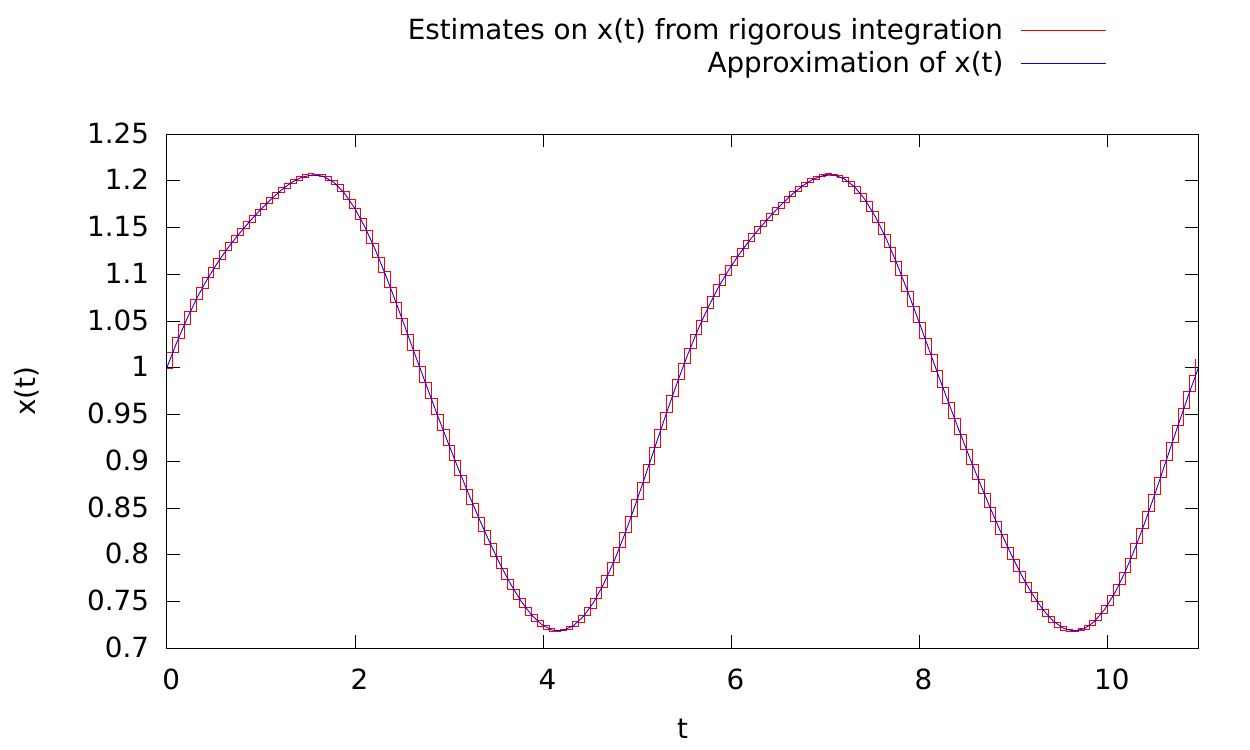}

\includegraphics[width=125mm]{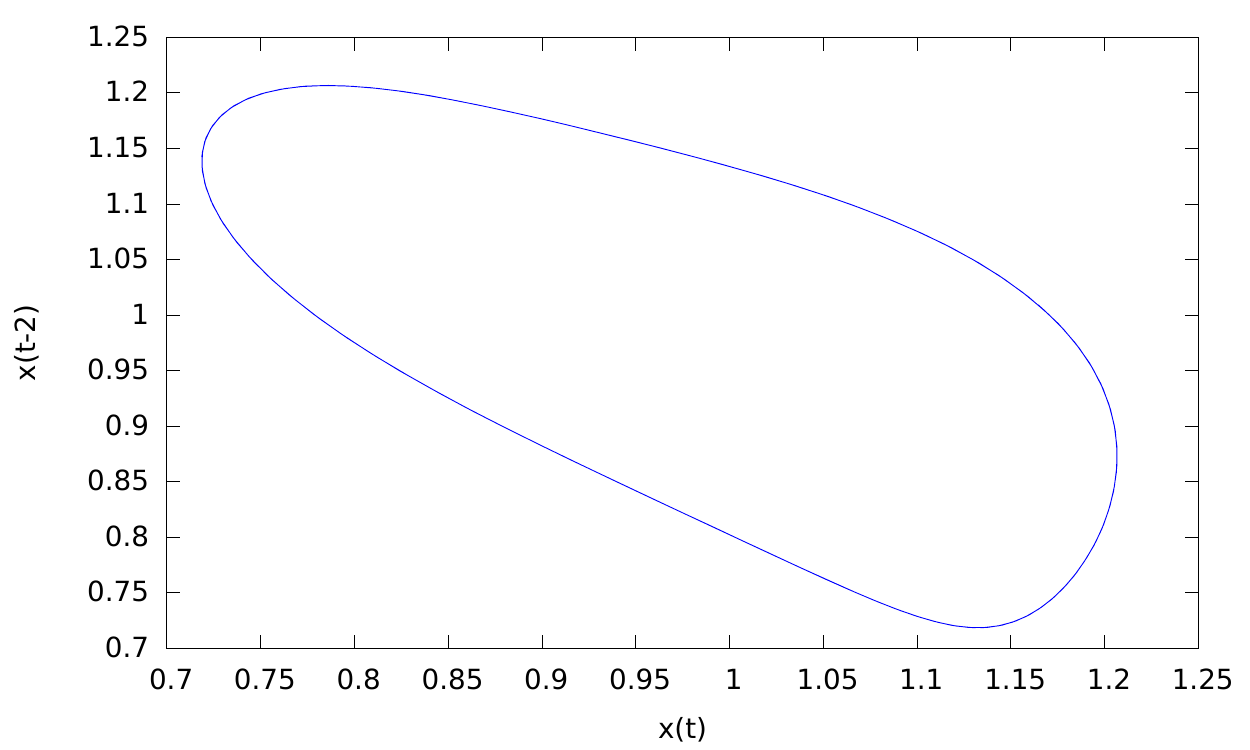}
}

\caption{\label{ fig:mackey-glass-n6-stable }Top: approximate function $\hat{x}$ (blue) and estimates on
the value of the true solution obtained from computer-assisted proof (red). Bottom: solution plotted
as parametric curve $r(t) = (\hat{x}(t), \hat{x}(t-\tau))$.
}.
\end{figure}

A graphical representation of the estimates obtained in the proof can be found in Figure~\ref{ fig:mackey-glass-n6-stable }.

\subsubsection{Case $n=8$}
For $n=8$ we consider the periodic orbit after the first period doubling. This time the period of
the orbit is long enough to overcome the initial loss of regularity, so  we consider the first
return Poincar\'e map.

Numerical computation shows that the orbit is attracting with the 10 most
significant eigenvalues of the map $P_{\ge \omega}$
estimated to be:

{\begin{center}
\begin{tabular}{|c|c|c|c|c|c|c|c|c|c|c|}
\hline
$Re \lambda$ & 0.3090 & -0.1359 & -0.0067 & $-7.58 \cdot 10^{-4}$ & $6.58 \cdot 10^{-4}$ & $-1.23 \cdot 10^{-4}$  & $ 2.184 \cdot 10^{-5}$ \\
\hline
$Im \lambda$ &        &         &         &                       &                      &                        & $ 6.265 \cdot 10^{-6}$ \\
\hline
$ |\lambda$| & 0.3090 &  0.1359 & 0.0067  & $ 7.58 \cdot 10^{-4}$ & $6.58 \cdot 10^{-4}$ & $ 1.23 \cdot 10^{-4}$  & $ 2.272 \cdot 10^{-5}$ \\
\hline
\end{tabular}
\end{center}
}

\begin{theorem}
\label{thm:mackey-glass-n8-stable}
There exists a $T$-periodic solution $x$ with period $T \in \left[11.1350, 11.1353\right]$
to Equation~\eqref{eq:mackey-glass} for parameters $\gamma = 2$, $\alpha = 1$, $\tau = 2$ and $n = 8$.
Moreover
\begin{eqnarray*}
    \left\| \hat{x} - x\right\|_{C^0} &\le& 0.012 \\
    \left\| \hat{x} - x\right\|_{C^1} &\le& 0.06 \\
    \left\| \hat{x} - x\right\|_{C^2} &\le& 0.20 \\
    \left\| \hat{x} - x\right\|_{C^3} &\le& 0.52 \\
    \left\| \hat{x} - x\right\|_{C^4} &\le& 1.25
\end{eqnarray*}
for $\hat{x}$ defined by
\begin{equation*}
\renewcommand*{\arraystretch}{1.5}    \begin{array}{rclclc}
\hat{x}(t) & = & 0.9480 & + & &  \\
           & + & 0.0477 \cdot \cos\left( {2\pi \over T} \cdot 1 \cdot t \right) & - & 0.0689 \cdot \sin\left( {2\pi \over T} \cdot 1 \cdot t \right) & + \\
           & + & 0.2516 \cdot \cos\left( {2\pi \over T} \cdot 2 \cdot t \right) & + & 0.1120 \cdot \sin\left( {2\pi \over T} \cdot 2 \cdot t \right) & + \\
           & + & 0.0242 \cdot \cos\left( {2\pi \over T} \cdot 3 \cdot t \right) & + & 0.0604 \cdot \sin\left( {2\pi \over T} \cdot 3 \cdot t \right) & - \\
           & - & 0.0386 \cdot \cos\left( {2\pi \over T} \cdot 4 \cdot t \right) & - & 0.0191 \cdot \sin\left( {2\pi \over T} \cdot 4 \cdot t \right) & + \\
           & + & 0.0132 \cdot \cos\left( {2\pi \over T} \cdot 5 \cdot t \right) & - & 0.0068 \cdot \sin\left( {2\pi \over T} \cdot 5 \cdot t \right) & - \\
           & - & 0.0197 \cdot \cos\left( {2\pi \over T} \cdot 6 \cdot t \right) & + & 0.0198 \cdot \sin\left( {2\pi \over T} \cdot 6 \cdot t \right) & + \\
           & + & 0.0077 \cdot \cos\left( {2\pi \over T} \cdot 7 \cdot t \right) & - & 0.0134 \cdot \sin\left( {2\pi \over T} \cdot 7 \cdot t \right) & - \\
           & - & 0.0053 \cdot \cos\left( {2\pi \over T} \cdot 8 \cdot t \right) & - & 0.0051 \cdot \sin\left( {2\pi \over T} \cdot 8 \cdot t \right) & - \\
           & - & 0.0005 \cdot \cos\left( {2\pi \over T} \cdot 9 \cdot t \right) & + & 0.0029 \cdot \sin\left( {2\pi \over T} \cdot 9 \cdot t \right) & - \\
           & - & 0.0018 \cdot \cos\left( {2\pi \over T} \cdot 10 \cdot t \right) & - & 0.0017 \cdot \sin\left( {2\pi \over T} \cdot 10 \cdot t \right) &
\end{array}
\end{equation*}
\end{theorem}
\textbf{Proof:} The proof follows the same lines as in the case of
Theorem~\ref{thm:mackey-glass-n6-stable} (except this time we consider the first return to the section).
Therefore we just list the parameters from the proof.
\begin{eqnarray*}
 p &=& 128, \quad n \ \ = \ \ 4 \\
    t_S                      &  \in & \omega + \left[\varepsilon\right]         \ \subset \ \left[11.1350,11.1353\right],    \\
    \omega                   &   =  & q \cdot \frac{\tau}{p}                    \     =   \ 712 \cdot \frac{\tau}{p},        \\
    \left[\varepsilon\right] &   =  & \left[0.01015698552,0.01016088515\right], \\
    l(\dot{x}) &\ge& 0.2636, \quad \mbox{for $x \in C^{n+1} \cap P(\bar{x}_0)$} \\
    \|\hat{x}-x\|_{C^0} &\leq& 0.01138319492 < 0.012.
\end{eqnarray*}
\qed

\begin{figure}
\centering{
\includegraphics[width=125mm]{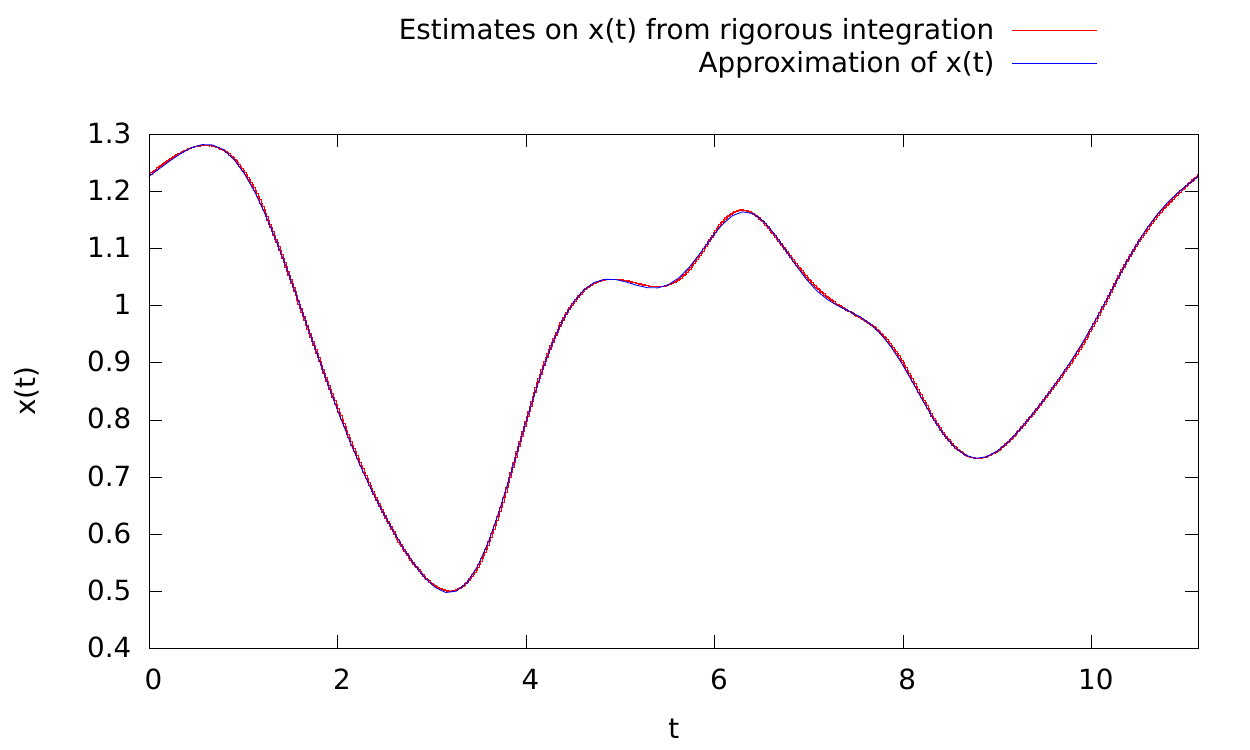}

\includegraphics[width=125mm]{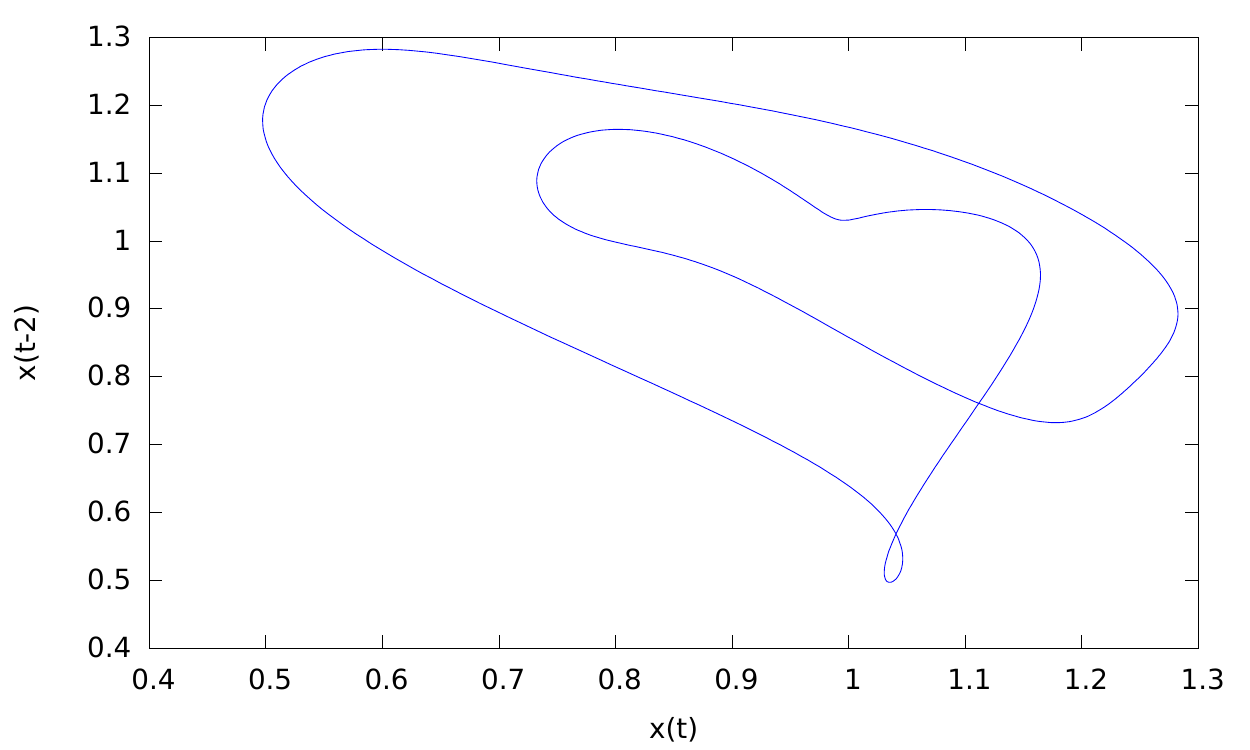}
}

\caption{\label{ fig:mackey-glass-n8-stable }Top: approximate function $\hat{x}$ (blue) and estimates on
the value of the true solution obtained from computer-assisted proof (red). Bottom: solution plotted
as parametric curve $r(t) = (\hat{x}(t), \hat{x}(t-\tau))$.
}.
\end{figure}

The diameter of the estimation for period $T$ (also for the last step $[\epsi_1, \epsi_2]$) obtained
from the computer-assisted proof is close to $3.899 \cdot 10^{-6}$.
A graphical representation of the estimates obtained in the
proof can be found in Figure~\ref{ fig:mackey-glass-n8-stable }.

The execution time was around $12$ minutes. This increase when compared to
$n=6$ is due to much larger representation size in this case which affects
the complexity of matrix and automatic differentiation algorithms
which we are using.

\section{Outlook and future directions}
\label{sec:conclusions}

The results presented in this work
might be improved in several ways:
\begin{itemize}
\item An extension of the integration algorithm to
the systems of delay equations in $\R^k$ for $k>1$.  This is rather straightforward and it does not require any new ideas;

\item A different representation of function
sets. Currently, we use the piecewise Taylor expansions, but other approaches,
like the Chebyshev polynomials, might  be better as they may produce
better approximations  on longer intervals;

\item avoiding the loss of
the regularity at the beginning of the integration, which imposes
the requirement for the transition time to section $t_S$ to be ,,long enough''.
The complete solution would be to confine the initial condition to the
invariant set $M^n \subset C^n$. We are currently working on this matter;
\end{itemize}

Other goal would be to apply the integrator to prove the existence
of hyperbolic periodic orbits with one or more unstable directions,
for example to establish the existence of LSOPs \cite{dde-lsop-krisztin}
in some general smooth DDEs, or unstable periodic solutions to Mackey-Glass equation.
Good theorems, suitable for that task, already exist, see \cite{zgliczynski-pde-3}
for the analogous question in the dissipative PDEs setting.

The ultimate goal is to establish tools to prove
chaotic dynamics in general DDEs, such as Mackey-Glass equation.

\section{Acknowledgements}
Research has been supported by Polish
National Science Centre grants 2011/03B/ST1/04780 and 2016/22/A/ST1/00077.

\bibliographystyle{plain} 	
\bibliography{dde}

\end{document}